\documentclass[11pt]{article}
\usepackage{amsmath}
\usepackage{latexsym}
\usepackage{amsfonts}
\usepackage{amssymb}
\usepackage{color}
\usepackage{amsmath, upgreek}
\usepackage{amssymb}
\usepackage{color}
\usepackage{amscd}
\usepackage{xspace}
\usepackage{verbatim}
\usepackage{graphicx}
\newcommand{\bth}[1]{\def\name{Theorem}
\begin{sub}\label{t:#1}}
\newcommand{\blemma}[1]{\def\name{Lemma}
\begin{sub}\label{l:#1}}
\newcommand{\bcor}[1]{\def\name{Corollary}
\begin{sub}\label{c:#1}}
\newcommand{\bdef}[1]{\def\name{Definition}
\begin{sub}\label{d:#1}}
\newcommand{\bprop}[1]{\def\name{Proposition}
\begin{sub}\label{p:#1}}

\newcommand{\BA}{\begin{array}}
\newcommand{\EA}{\end{array}}
\newcommand{\BAN}{\renewcommand{\arraystretch}{1.2}
\setlength{\arraycolsep}{2pt}\begin{array}}
\newcommand{\BAV}[2]{\renewcommand{\arraystretch}{#1}
\setlength{\arraycolsep}{#2}\begin{array}}
     \def\gb{\beta}       \def\gg{\gamma}
       \def\gd{\delta}      \def\ge{\epsilon}
\def\gth{\theta}                         
\def\gf{\phi}           
      \def\gk{\kappa}      \def\gl{\lambda}
\def\gm{\mu}        \def\gn{\nu}         \def\gp{\pi}
            
\def\gs{\sigma}       \def\gt{\tau}
   \def\gv{\vartheta}   \def\gw{\omega}
                \def\gz{\zeta}
\def\Gg{\Gamma}     \def\Gd{\Delta}      \def\Gf{\Phi}
\def\Gth{\Theta}
\def\Gl{\Lambda}    \def\Gs{\Sigma}      
\def\Gw{\Omega}              
   \def\CM{{\mathcal M}}   
   \def\CO{{\mathcal O}}   
      
\def\CD{{\mathcal D}}      
      
      \def\CL{{\mathcal L}}

\def\CW{{\mathcal W}} 

   \def\BBN {\mathbb N}    
   \def\BBR {\mathbb R}

\newcommand{\R}{\mathbb{R}}

\newcommand{\N}{\mathbb{N}}

\newcommand{\beq}{\begin{equation} }
\newcommand{\eqq}{\end{equation} }
\newcommand{\cuad}{{\sqcap\kern-.68em\sqcup}}
\newcommand{\abs}[1]{\mid #1 \mid}
\newcommand{\norm}[1]{\|#1\|}

\newtheorem{definition}{Definition}[section]
\newtheorem{teo}{Theorem}[section]

\newtheorem{proposition}{Proposition}[section]

\newtheorem{lemma}{Lemma}[section]
\newtheorem{corollary}{Corollary}[section]
\newtheorem{remark}{Remark}[section]
\newcommand{\bremark}{\begin{remark} \em}
\newcommand{\eremark}{\end{remark} }

\newcommand{\qeda}{\hspace{10mm}\hfill $\square$}

\def\beeq{\begin{equation}}
\def\eeq{\end{equation}}
\newcommand{\begeqaet}{\begin{eqnarray*}}
\newcommand{\eneqaet}{\end{eqnarray*}}

\begin{document}

\begin{center}{\bf  \large
Qualitative properties of  solutions to semilinear \\[2mm]

  elliptic equations  from the gravitational \\[2mm]

 Maxwell  Gauged $O(3)$ Sigma model  }
\bigskip
\bigskip
{\small

  {\sc Huyuan Chen\footnote{chenhuyuan@yeah.net},\quad  Hichem Hajaiej \footnote{hhajaie@calstatela.edu}\quad and\quad Laurent V\'eron\footnote{ Laurent.Veron@lmpt.univ-tours.fr} }
  \bigskip

  $^1$ Department of Mathematics, Jiangxi Normal University,\\
  Nanchang, Jiangxi 330022, PR China \\[12pt]

$^2$ California State University, Los Angeles,\\
 5151 University Drive, Los Angeles, CA 90032-8530, USA\\[12pt]

$^3$ Laboratoire de Math\'{e}matiques et Physique Th\'{e}orique Universit\'{e} de Tours, 
\\ 37200 Tours, France
}
  \bigskip

\begin{abstract} 
This article is devoted to the study of the following semilinear equation with measure data which originates in  the gravitational
 Maxwell  gauged $O(3)$ sigma model,
$$
-\Delta u +A_0(\prod^k_{j=1}|x-p_j|^{2n_j} )^{-a} \frac{e^u}{(1+e^u)^{1+a}}=4\pi\sum_{j=1}^k n_j\delta_{p_j}-4\pi\sum^l_{j=1}m_j\delta_{q_j} \quad{\rm in}\;\; \mathbb{R}^2.\leqno(E)
$$
In this equation the  $\{\delta_{p_j}\}_{j=1}^k$ (resp. $\{\delta_{q_j}\}_{j=1}^l$ ) are Dirac masses  concentrated at the points
 $\{p_j\}_{j=1}^k$, (resp. $\{q_j\}_{j=1}^l$),
$n_j$ and $m_j$ are positive integers, and $a$ is a nonnegative real number. We set  $ N=\sum^k_{j=1}n_j $ and $  M= \sum^l_{j=1}m_j$.\smallskip

In previous works \cite{C,Y2}, some qualitative properties of solutions of $(E)$ with $a=0$ have been established.  Our aim in this article is to study  the more general  case where $a>0$.  The additional difficulties of this case come from the fact that the nonlinearity is no longer monotone and the data are signed measures. As a consequence  we cannot anymore construct directly the solutions by the monotonicity method combined with the  supersolutions and subsolutions technique. Instead we develop a new and self-contained approach which enables us to emphasize the role played by the gravitation in the gauged $O(3)$ sigma model.
Without the gravitational term, i.e. if $a=0$, problem $(E)$ has a layer's structure of solutions $\{u_\beta\}_{\beta\in(-2(N-M),\, -2]}$, where  $u_\beta$ is the unique non-topological solution such that $u_{\beta}=\beta\ln |x|+O(1)$ for $-2(N-M)<\beta<-2$ and $u_{-2}=-2\ln |x|-2\ln\ln |x|+O(1)$ at infinity respectively.
On the contrary,  when $a>0$,     the set of solutions to problem $(E)$ has a much richer structure: besides the topological solutions, there exists 
 a sequence of non-topological solutions in type I, i.e. such that $u $ tends to $-\infty$  at infinity, and of non-topological solutions of type II,
which tend to $\infty$ at infinity. The existence of these types of solutions depends on the values of the parameters $N,\, M,\, \beta$ and on the gravitational interaction associated to $a$. 

\end{abstract}\end{center}

\vspace{0.7mm}
  \noindent {\bf Keywords}:    Gauged Sigma Model;  Non-topological Solution;  Topological Solution.

  \smallskip

\noindent {\small {\bf MSC2010}: 35R06, 35A01, 	81T13. }

\tableofcontents
\hspace{.05in}

\vspace{2mm}

\setcounter{equation}{0}
\section{Introduction}
In this paper our goal  is to  classify  the solutions of the following equation with measure data
\begin{equation}\label{eq 1.1}
-\Delta u +A_0(\prod^k_{j=1}|x-p_j|^{2n_j} )^{-a} \frac{e^u}{(1+e^u)^{1+a}} =4\pi\sum_{j=1}^k n_j\delta_{p_j}-4\pi\sum^l_{j=1}m_j\delta_{q_j} \quad{\rm in}\;\;  \mathbb{R}^2,
\end{equation}
where  $\{\delta_{p_j}\}_{j=1}^k$ (resp. $\{\delta_{q_j}\}_{j=1}^l$ ) are Dirac masses  concentrated at the points
 $\{p_j\}_{j=1}^k$, (resp. $\{q_j\}_{j=1}^l$), $p_j\not=p_{j'}$ for $j\not=j'$, the related coefficients $n_j$ and $m_j$ are positive integers, $A_0>0$ is a given constant, $a=16\pi G$ with $G$ being the Newton's gravitational constant (or more precisely a dimensionless
rescaling factor of the gravitational constant \cite{Y1}) which is of the order of $10^{-30}$. This means that physically speaking the exponent $a$ is very small. Set
\begin{equation}\label{eq 1.1x}
{\bf P}(x)=A_0(\prod^k_{j=1}|x-p_j|^{2n_j} )^{-a}.
\end{equation}
Since 
\begin{equation}\label{eq 1.1y}
2^{-1-a}\min\{e^u,e^{-au}\}\leq \frac{e^u}{(1+e^u)^{1+a}}\leq \min\{e^u,e^{-au}\},
\end{equation}
we define the notion of weak solution as follows:
\begin{definition}\label{ws} A function $u\in L^1_{loc}(\R^2)$ such that  ${\bf P}\min\{e^u,e^{-au}\}\in L^1_{loc}(\R^2)$ is called  a weak  solution of (E), if for any $\xi\in C_c^\infty(\R^2)$,
$$\int_{\R^2} u(-\Delta)\xi\, dx+\int_{\R^2} {\bf P} \frac{e^u}{(1+e^u)^{1+a}} \xi\, dx=4\pi\sum^{k}_{j=1}n_j\xi(p_j)-4\pi\sum^{l}_{j=1}m_j\xi(q_j). $$
\end{definition}
This definition means that the following equation holds in the sense of distributions in $\mathbb{R}^2$,
\begin{equation}\label{eq 1.1a}
-\Delta u +{\bf P}\frac{e^u}{(1+e^u)^{1+a}} =4\pi\sum_{j=1}^k n_j\delta_{p_j}-4\pi\sum^l_{j=1}m_j\delta_{q_j}.
\end{equation}

We denote by $\Sigma:=\{p_1,\cdots,p_k,q_1,\cdots, q_l\}$ the set of the supports of the measures in the right-hand side of (\ref{eq 1.1}).  
Since the nonlinearity  in (\ref{eq 1.1a}) is locally bounded in $\R^2\setminus \Sigma$,  a weak solution of (\ref{eq 1.1a}) belongs to $C^2(\R^2\setminus \Sigma)$ and is a strong solution of 
\begin{equation}\label{eq 1.1b}
-\Delta u +{\bf P} \frac{e^u}{(1+e^u)^{1+a}} =0\quad\mbox{in }\;\R^2\setminus \Sigma.
\end{equation}
The nonlinear term is not monotone, actually the function $u\mapsto \frac{e^u}{(1+e^u)^{1+a}}$ is increasing on $(-\infty, -\ln a)$, and decreasing on 
$( -\ln a,\infty)$. This makes the structure of solutions of our problem much more complicated  than the case where $a=0$. 
\subsection{Physical models and related equations}

Equation (\ref{eq 1.1}) comes from the Maxwell gauged $O(3)$ sigma model.
 When $a=0$, it governs the self-dual $O(3)$ gauged sigma model  developed from  Heisenberg ferromagnet, see references \cite{B,BP,R,S}.   
When the sigma model for Heisenberg ferromagnet with magnetic field is two-dimensional, it can be expressed by a local $U(1)$-invariant action density \cite[p. 43-49]{Y2}:
$$\mathcal{L} =-\frac14 F_{\mu\nu}F^{\mu\nu}+\frac12 D_\mu \phi \overline{D^\mu \phi}  -\frac12(1-\vec{n} \cdot \phi)^2,$$
where $\vec{n}=(0,0,1)$, $\phi=(\phi_1,\phi_2,\phi_3)$  is a spin vector  defined over the $(2+1)$-dimensional Minkowski spacetime $\R^{2,1}$, with value in the unit sphere $ \mathbb{S}^2$, i.e. $|\phi|=1$, $D_\mu$ are gauge-covariant derivatives on $\phi$, defined by
$$D_\mu \phi =\partial _\mu \phi +A_\mu(\vec{n} \times \phi) \quad \text{where }\;\mu=0,1,2$$
and $F_{\mu\nu}=\partial_\mu A_\nu-\partial_\nu A_\mu$ is the electromagnetic curvature induced from
the 3-vector connection $A_\nu$, $\gn=0,1,2$ as detailled in \cite[p. 177-189]{Y4}. When 
the time gauge $A_0$ is zero, that is in the static situation, the functional of total energy can be expressed by the following expressions
\begin{eqnarray*}
E(\phi,A) &=& \frac12\int_{\R^2}\left((D_1\phi)^2+(D_2\phi)^2+ (1-\vec{n} \cdot \phi)^2+F_{12}^2\right)dx \\
   &=& 4\pi|deg(\phi)|+\frac12 \int_{\R^2}\left((D_1\phi\pm \phi\times D_2\phi)^2+(F_{12}\mp (1-\vec{n} \cdot \phi))^2\right)dx,
\end{eqnarray*}
where $deg(\phi)$ denotes the Brouwer's degree of $\phi$. The related Bogomol'nyi equation is obtained  by using  the stereographic projection $\phi\mapsto\tilde\phi$ from
the south pole $\mathcal{S}=(0,0,-1)$   of $\mathbb{S}^2\setminus \{\mathcal{S}\}$ onto $\R^2$ (see e.g. \cite{Chae,Y4} for details). Then the function $u=\ln|\tilde\phi|^2$ satisfies
\begin{equation}\label{eq 1.2}
-\Delta u +  \frac{4e^u}{1+e^u } =4\pi\sum_{j=1}^k n_j\delta_{p_j}-4\pi\sum^l_{j=1}m_j\delta_{q_j} \quad{\rm in}\;\;  \mathbb{R}^2.
\end{equation}
It is  pointed out in \cite{Y2}   that the points $p_j$ ($j=1,\cdots, k$), which are the  poles of $\tilde\phi$
can be viewed as magnetic monopoles and the points $q_j$ ($j=1,\cdots, l$), which are the zeros of $\tilde\phi$   as  antimonopoles (see \cite[p. 55]{Y4}).
They are also called  magnetic vortices and anti-vortices respectively.

An important quantity for the gauged sigma model is the total magnetic flux. It is customary \cite{Sc} to identity it to the integral of the curvature as follows:
\begin{equation}\label{tmf-0}
 \mathcal{M}(\phi)=\int_{\R^2}F_{12}.
\end{equation}
Using the variable $u$ its value coincides with $\int_{\R^2}\Gd udx$ (the Laplacian being taken a.e.). Thus, for the sake of simplicity, we identify  $\mathcal{M}(\phi)$ and  $\mathcal{M}(u)$, an expression which will be called the total flux in the sequel. 
Here and in what follows, we denote
$$N=\sum^k_{j=1}n_j\quad{\rm and}\quad  M= \sum^l_{j=1}m_j.$$

When the gravitation constant $G$ is replaced by zero, a layer's structure of solutions of (\ref{eq 1.1}) 
has been determined in the following result:
\begin{teo}\label{teo 3.1} \cite{C,Y2}
$(i)$ If $M=N-1$, then problem (\ref{eq 1.2}) has no solution.\smallskip

\noindent$(ii)$ If $M<N-1$,  then for any $\beta\in[2,\,2(N-M))$ problem (\ref{eq 1.2}) has a unique solution $u_\beta$ verifying
$$\mathcal{M}(u_\beta) =2\pi(2(N-M)+\beta),$$
with the following behaviour as $|x|\to\infty$,
$$u_\beta(x) =\arraycolsep=1pt\left\{
\begin{array}{lll}
 \displaystyle -\beta\ln|x|+O(1)\quad&{\rm if} \ \    \beta\in(2,2(N-M)),
\\[2mm]\phantom{ }
 \displaystyle  -2\ln|x|-2\ln\ln|x|+O(1)\quad&{\rm if} \ \    \beta=2.
\end{array}
\right.$$
Furthermore the correspondence $\gb\mapsto u_\gb$ is decreasing.\smallskip

\noindent$(iii)$ If  $M<N-1$ and $u$ is a non-topological solution of (\ref{eq 1.2}) with finite total magnetic flux, i.e. $\mathcal{M}(u)<\infty$,
then  there exists a unique $\beta\in[2,2(N-M))$ such that $u=u_\beta$.

\end{teo}

These equations have been studied extensively, motivated by a large range of many applications in physics  such as the gauged sigma models with broken symmetry \cite{Y3}, the gravitational Maxwell gauged $O(3)$ sigma model \cite{Chae,CY,Sc,S},  the self-dual Chern-Simons-Higgs model \cite{CFL,LPY}, magnetic vortices \cite{LY}, Toda system \cite{LWY,PT}, Liouville equation \cite{JW} and the references therein. It is also motivated by important questions in the theory of nonlinear partial differential equations \cite{BM,Va,V}, which has its own
features in two dimensional space.


When $a=16\pi G$, equation (\ref{eq 1.1}) governs the gravitational Maxwell  gauged $O(3)$ sigma model  restricted to a plane. Because of the gravitational interaction between particles, the Lagrangian density becomes
$$\mathcal{L} = \frac14g^{\mu \mu'}g^{\nu\nu'} F_{\mu\nu}F_{\mu'\nu'}+\frac12 D_\mu \phi \overline{D^\mu \phi}  -\frac12(1-\vec{n} \cdot \phi)^2$$
with stress energy tensor
$$T_{\mu\nu}=g^{\mu'\nu'} F_{\mu\nu}F_{\mu'\nu'}+D_\mu \phi  D_\mu \phi-g_{\mu\nu}\mathcal{L}. $$
We simplify the Einstein equation
$$R_{\mu\nu}-\frac12Rg_{\mu\nu}=-8\pi GT_{\mu\nu},$$
where $R_{\mu\nu}$ is the Ricci tensor and $R$ is a scalar tensor of the metric in considering a metric conformal to the (2 + 1)-dimensional Minkowski one
 $$g_{\mu\nu}= \left(
                     \begin{array}{ccc}
                       -1 & 0 & 0 \\
                       0 & e^\eta &  0\\
                       0 & 0& e^\eta \\
                     \end{array}
                   \right).$$
Then
$$\frac12 e^{-\eta} \Delta \eta=-8\pi GT_{00},$$
where
$$T_{00}=\frac12(e^{-\eta}F_{12}\pm(1-\vec{n}\cdot \phi))^2\pm e^{-\eta} F_{12} (1-\vec{n}\cdot \phi)) \pm e^{-\eta} \phi(D_1\phi\times D_2\phi)+\frac12(D_1\phi\pm \phi\times D_2\phi)^2. $$
The minimum of the energy is achieved if and only if $(\phi,A)$ satisfies the self-dual equations (the Bogomol'nyi equations)
$$D_1\phi=\mp \phi\times D_2\phi,\quad F_{12}=\pm e^\eta (1-\vec{n}\cdot\phi).$$
Furthermore, a standard analysis   yields equation (\ref{eq 1.1}). In particular, Yang in \cite{Y4} studied  equation (\ref{eq 1.1}) when there is only one concentrated   pole, i.e. $k=1$ and $l=0$. For multiple poles, Chae  showed in \cite{Chae} that problem (\ref{eq 1.1}) has a sequence of non-topological solutions $u_\beta$
such that
$$u_\beta(x)=\beta\ln|x|+O(1)\quad{\rm when }\; |x|\to\infty$$
for $\beta\in(-\min\{6,\, 2(N-M)\},\, -2)$, when
\begin{equation}\label{assump 1}
 aN<1 \quad{\rm and}\quad N-M\geq 2.
\end{equation}
Under the assumption (\ref{assump 1}), the existence of solutions has been improved up to the range $\beta\in(- 2(N-M),\, -2)$ by Song
in \cite{S}. However, these existence results do not show the role of the gravitation played in the gauged sigma model and the features of the interaction of the diffusion and the non-monotone nonlinearity of equation (\ref{eq 1.1}) in the whole two dimensional space.

\subsection{Main results}

Note that if we take into account the gravitation,  the total magnetic flux turns out to be
\begin{equation}\label{tmf}
  \mathcal{M}(u) =\int_{\R^2} {\bf P}(x)\frac{e^u}{(1+e^u)^{1+a}}  \, dx,
 \end{equation}
which, due to the potential and the decay to zero for
$\frac{e^t}{(1+e^t)^{1+a}}$ as $t\to \infty$, allows the existence of solutions with very wild behaviors at infinity.  In fact,  the following three  types of solutions are considered
in this paper
$$
\arraycolsep=1pt\left\{
\begin{array}{lll}
\mbox{  a solution $u$ of (\ref{eq 1.1}) is  topological   } \quad &{\rm if } \ \ \displaystyle \lim_{|x|\to+\infty} u(x)=\ell\in\R,\\[1mm]
 \phantom{ }
 \mbox{  a solution $u$ of (\ref{eq 1.1}) is  non-topological of type I } \quad &{\rm if } \ \  \displaystyle\lim_{|x|\to+\infty} u(x)=-\infty,\\[1mm]
 \phantom{ }
  \mbox{  a solution $u$ of (\ref{eq 1.1}) is  non-topological of type II } \quad &{\rm if } \ \  \displaystyle\lim_{|x|\to+\infty} u(x)=+\infty.
\end{array}
\right.$$
\smallskip

  The first  result of this paper deals with non-topological solutions of type I for (\ref{eq 1.1}). For such a task we introduce two important quantities:
 \begin{equation}\label{1.1+q}
\beta^\#= \max\Big\{-2(N-M),\,\frac{2-2aN}{a}\Big\}\quad{\rm and}\quad   \beta^*=\min\left\{0,\, 2aN-2, \alpha^*-2(N-M) \right\},
\end{equation} 
where
\begin{equation}\label{1.1+p}
 \alpha^*:= \frac{1}{2\pi} \int_{\R^2}{\bf P}(x) dx.
\end{equation}

Notice that  $\alpha^*=\infty$ if  $a n_j\geq 1$ for some $j$ or 
$aN \le 1$,
 otherwise $\alpha^*$ is finite, in this case, a free parameter $A_0$ should be taken into account. If $aN\le 1$, we have that  $\beta^* =2aN-2\leq 0$.

\begin{teo}\label{teo 0}
Let $a=16\pi G$, $an_j< 1\ \, {\rm for}\ j=1,\cdots,k$ and $\mathcal{M}$ be the total magnetic flux  given in (\ref{tmf}).\smallskip

\noindent $(i)$ If
\begin{equation}\label{NM1}
  aN\le 1 \quad{\rm and}\quad  M<(1+a)N-1,
\end{equation}
  then   for any $\beta\in(-2(N-M),\, \beta^*)$,  problem (\ref{eq 1.1}) possesses a minimal  solution $u_{\beta,min}$
  satisfying
$$u_{\beta,min}(x)=\beta\ln |x|+O(1) \quad{\rm as}\quad |x|\to+\infty.$$
Moreover, for some real number $C_*$,
\begin{equation}\label{1.3}
  u_{\beta,min}(x)=\beta\ln|x|+C_*+O(|x|^{-\frac{aN-\beta-2}{aN-\beta-1}})\quad{\rm as}\quad |x|\to+\infty,
\end{equation}
and the  total magnetic flux of the solution  $ u_{\beta,min} $ is equal to $2\pi[2(N-M)+\beta]$, i.e.
  \begin{equation}\label{1.2}
   \mathcal{M}(u_{\beta,min}) =2\pi[2(N-M)+\beta].
  \end{equation}

\noindent$(ii)$ If
\begin{equation}\label{1.5}
 aN>1 \quad{\rm and}\quad  M<N,
\end{equation}
then  $\beta^\# <0$ and for any $ \beta\in\left(\beta^\#  ,\, 0\right)$,
  problem (\ref{eq 1.1}) possesses a sequence of non-topological  solutions $u_{\beta,i}$ of type I satisfying
\begin{equation}\label{1.7-1}
  u_{\beta,i}(x)=\beta\ln|x|+C_i+O(|x|^{-\frac{2aN-2\beta-2}{2aN-2\beta-1}})\quad{\rm as}\quad |x|\to\infty,
\end{equation}
where
$$C_i<C_{i+1}\to\infty\quad{\rm as}\quad i\to+\infty.$$
Moreover, the  total magnetic flux of the solutions $\{u_{\beta,i}\}_i$ is equal to $2\pi[2(N-M)+\beta]$.
 \end{teo}

Note that our assumption (\ref{NM1}) is much weaker than (\ref{assump 1}) and Theorem \ref{teo 0} provides a  larger range of $\beta$ for existence of
solutions $u_\beta$ verifying $u_\beta=\beta\ln |x|+o(1)$ at infinity. Furthermore we obtain a minimal solution and not just a finite energy solution as in 
\cite[Theorem 1.3]{S}. Note also that the assumption $M<(1+a)N-1$ implies that $ \beta^*>-2(N-M)$,
and our second interest is to consider this extremal case $\beta=\beta^*$, which 
is $2aN-2$ under the assumption (\ref{NM1}).

 \begin{teo}\label{teo 0-c}
 Assume that $a=16\pi G$,   $an_j< 1\ \, {\rm for}\ j=1,\cdots,k$, the magnetic flux $ \mathcal{M}$ is given by (\ref{tmf}) and  let (\ref{NM1}) hold.
 
 \noindent Then    problem (\ref{eq 1.1}) possesses a minimal non-topological solution $u_{\beta^*,min}$
   satisfying
   \begin{equation}\label{1.3-1}
  u_{\beta^*,min}(x)=\beta^*\ln|x|-2\ln\ln|x|+O(1)\quad {\rm as} \;    |x|\to+\infty, 
\end{equation}
and  the  total magnetic flux of $u_{\beta^*,min}$   is equal to $2\pi[2(N-M)+\beta^*]$.
 \end{teo}

The existence of non-topological states of type II to (\ref{eq 1.1}) states as follows.

 \begin{teo}\label{teo 1}
Assume that  $a=16\pi G$,  $an_j< 1\ \, {\rm for}\ j=1,\cdots,k$  and $\beta^\#$ is given by (\ref{1.1+q}),
 then  ${\rm for\ any \ } \ \beta>  \beta^\#_+=\max\{0,\,\beta^\#\},$
  problem (\ref{eq 1.1}) possesses a sequence of non-topological  solutions $\{u_{\beta,i}\}_i$ such that
\begin{equation}\label{1.7}
  u_{\beta,i}(x)=\beta\ln|x|+C_i+O(|x|^{-\frac{2aN-2\beta-2}{2aN-2\beta-1}})\quad {\rm as} \;    |x|\to+\infty,
\end{equation}
where
$$C_i<C_{i+1}\to+\infty\quad {\rm as} \;    i\to+\infty.$$
Moreover, the  total magnetic flux of the solutions $\{u_{\beta,i}\}_i$ is equal to $2\pi[2(N-M)+\beta]$.


 \end{teo}

Concerning topological solutions  of (\ref{eq 1.1}), we have following result, 

 \begin{teo}\label{teo 2}
Let $a=16\pi G$,  $an_j< 1\ \, {\rm for}\ j=1,\cdots,k$ and (\ref{1.5}) hold true.
 
 \noindent Then   problem (\ref{eq 1.1}) possesses infinitely many topological solutions $u_{0,i}$ satisfying 
\begin{equation}\label{1.6}
 u_{0,i}(x) = C_i+O(|x|^{-\frac{2aN-2}{2aN-1}})\quad {\rm as} \;    |x|\to\infty,
\end{equation}
 where
 $$C_i<C_{i+1}\to\infty\quad {\rm as} \    i\to\infty.$$
Moreover, the  total magnetic flux of the solutions $\{u_{0,i}\}_i$ is equal to $4\pi(N-M)$.
\end{teo}

Note that Theorem \ref{teo 1} and Theorem \ref{teo 2} provide respectively infinitely many non-topological solutions of Type II and topological solutions. Furthermore,
there is no upper bound for these solutions, this is  due to the failure of the Keller-Osserman condition for the nonlinearity $\frac{4e^u}{(1+e^u)^{1+a}}$, see \cite{K,O}. More precisely equation (\ref{eq 1.1}) admits no solution with boundary blow-up in a bounded domain.   The existence of these solutions
illustrates that the gravitation plays an important role in the
 Maxwell  gauged $O(3)$ sigma model: \smallskip
 
 \noindent $(i)$ the set of solutions is extended to topological and two types of non-topological solutions; \smallskip
 
  \noindent $(ii)$ the uniqueness fails for the solution under the given condition $u_\beta(x)=\beta\ln |x|+O(1)$ at infinity; \smallskip
  
   \noindent $(iii)$ the numbers (counted with multiplicity) of magnetic poles $N,M$ do no longer verify $M<N+1$. In fact, for the non-topological solution of type I, it becomes $M<(1+a)N+1$, but   for the non-topological solution of type II, there is no restriction on $N$ and $M$, if $\beta>0$ is large enough. \smallskip
   
   \noindent Our existence statements of solutions of (\ref{eq 1.1}) are summarized in the three tables above.\smallskip

{\small
\begin{table}[tph]
\centering  
\caption{ Non-topological solutions of Type I  }\smallskip
\begin{tabular}{|c|c|c|c|}  
\hline
   Assumptions on $a, N, M$           &  range of $\beta$                              & solutions      & asymptotic behavior at $\infty$         \\ \hline
  $aN\le 1$,  $M<(1+a)N-1$         &  $(-2(N-M),\, \beta^*)$   & Minimal        & $\beta\ln|x| +O(1)$                    \\ \hline
  $aN\le 1$, $M<(1+a)N-1$            &   $ \beta^*=2(aN-1)$                &  Minimal       &  $\beta^*\ln|x| -2\ln\ln |x|+O(1)$    \\ \hline  
 $N>M$,  $aN> 1$           &  $ ( \beta^\#,\, 0 )$   & Multiple        & $\beta\ln|x|+c_i +o(1)$, $\displaystyle\lim_{i\to\infty}c_i=\infty$                    \\ \hline
\end{tabular}\medskip
\caption{ Non-topological solutions of Type II  }\smallskip
\begin{tabular}{|c|c|c|}  
\hline
  range of $\beta$                                                              & solutions      & asymptotic behavior at $\infty$         \\ \hline
 $( \beta^\#_+,\,\infty)$           & Multiple        & $\beta\ln|x|+c_i +o(1)$, $\displaystyle\lim_{i\to\infty}c_i=\infty$     \\ \hline
\end{tabular}\medskip
\caption{ Topological solutions  }\smallskip
\begin{tabular}{|c|c|c|}  
\hline
 Assumptions on $a, N, M$                      & solutions            & asymptotic behavior at $\infty$         \\ \hline
$aN>1$, $M<N$              &  Multiple       &  $ c_i+o(1)$, $\displaystyle\lim_{i\to\infty}c_i=\infty$
\\ \hline
\end{tabular}
\end{table}}

The biggest difference with the case that $a=0$ is that the nonlinearity is no longer monotone,  which makes more difficult to construct super and sub solutions to (\ref{eq 1.1}). Our main idea is to approximate
the solution by monotone iterative schemes for some related equations with an increasing nonlinearity.

Finally, we concentrate on the nonexistence of solutions  $u_\beta$ for (\ref{eq 1.1})
with the behavior $\beta\ln |x|+O(1)$ at infinity for some $\beta$.

 \begin{teo}\label{teo 3}
Assume that  $a=16\pi G$ and  $an_j< 1\ \, {\rm for}\ j=1,\cdots,k$. \smallskip

\noindent $(i)$ If  $aN<1$ and $\beta^*<\beta< \frac{2-aN}{a}$,
then problem (\ref{eq 1.1}) has no solution $u_\beta$ with the aymptotic behavior 
 $$u_\beta(x) =\beta\ln|x|+o(\ln|x|)\quad {\rm as} \ \    |x|\to\infty.$$

\noindent $(ii)$ If  $aN=1,$
then   problem (\ref{eq 1.1}) has no topological solution.


 \end{teo}

The remaining of this paper is organized as follows.  In Section 2, we present some decompositions of solutions of (\ref{eq 1.1}), some important estimates are provided and related forms of equations  are considered. We prove  that problem \ref{eq 1.1} has a minimal non-topological solution of Type I and minimal solutions in Section 3. Existence of infinitely many non-topological solutions of Type II is obtained in Section 4. Infinitely many topological solutions and minimal topological solution are constructed in Section 5.
 Finally, Section 6 deals with the classification of general non-topological solutions of (\ref{eq 1.1}) with infinite total magnetic flux.

\setcounter{equation}{0}
\section{Preliminary }

\subsection{Regularity}
We begin our analysis by  considering the regularity of  weak solutions of (\ref{eq 1.1}).
 Let $\zeta$ be a smooth and increasing function defined in $(0,\infty)$ and such that
$$\zeta(t)=\arraycolsep=1pt\left\{
\begin{array}{lll}
\ln t\quad &{\rm for } \ \ 0<t\le 1/2,
\\[2mm]\phantom{ }
0\quad  &{\rm for } \ \  t\geq 1.
\end{array}
\right.$$
Set
\begin{equation}\label{nu}
 \nu_1(x)=2\sum_{i=1}^kn_i\zeta\left(\tfrac{|x-p_i|}{\sigma}\right)\quad{\rm and}\quad  \nu_2(x)=2\sum_{j=1}^lm_j\zeta\left(\tfrac{|x-q_j|}{\sigma}\right),
\end{equation}
where $\sigma\in(0,1)$ is chosen such that any two balls of the set
 $$\left\{B_{\sigma}(p_i),\, B_{\sigma}(q_j):\, i=1,\cdots k,\,   j=1,\cdots l\right\}$$ 
 do not intersect.
We fix  a positive number $r_0\geq e^e$  large enough such that $B_{\sigma}(p_i), B_{\sigma}(q_j)\subset B_{r_0}(0)$ for $i=1,\cdots, k$ and $j=1,\cdots, l$, and we denote
$$\Sigma_1=\{p_1,\cdots,p_k\},\quad  \Sigma_2=\{q_1,\cdots,q_l\}\quad{\rm and}\quad \Sigma=\Sigma_1\cup \Sigma_2.$$
If $u$ is a weak solution of (\ref{eq 1.1}), we set
\begin{equation}\label{3.1}
u= w-\nu_1+\nu_2\quad{\rm in}\;\; \R^2\setminus \Gs.
\end{equation}
Then $w$ is a weak solution of
\begin{equation}\label{eq 3.2}
\displaystyle - \Delta w +{ V}\frac{  e^{w}}{(e^{\nu_1-\nu_2}+  e^{w})^{1+a}} =f_1-f_2\quad{\rm in}\;\; \R^2,
\end{equation}
with
\begin{equation}\label{V}
{ V}={\bf P}e^{a(\gn_1-\gn_2)}\,,\; f_1= 4\pi \displaystyle\sum_{i=1}^kn_i\delta_{p_i}-\Delta \nu_1\;\text{ and }\;\, f_2= 4\pi\displaystyle\sum_{j=1}^l m_j\delta_{q_j}-\Delta \nu_2.
\end{equation}
The functions $f_1,f_2$  are  smooth with compact supports in $B_{r_0}(0)$ and they satisfy
\begin{equation}\label{Vf}\int_{\R^2}(f_1-f_2) \, dx=4\pi(N-M).\end{equation}

\begin{proposition}\label{pr 2.1}
Assume that $u$ is a weak solution of (\ref{eq 1.1}), then $u$ is a classical solution of
\begin{equation}\label{eq 2.2}
-\Delta u +{\bf P}(x)\frac{e^u}{(1+e^u)^{1+a}} =0 \quad{\rm in}\;\;  \mathbb{R}^2\setminus\Sigma,
\end{equation}
and $w=u-\gn_1+\gn_2$ is a classical solution of (\ref{eq 3.2}) in whole $\mathbb{R}^2$.
\end{proposition}

\noindent{\bf Proof.} Let $u$ be a weak solution of (\ref{eq 1.1}). Since $\frac{e^u}{(1+e^u)^{1+a}}$ is uniformly bounded in $\mathbb{R}^2$ and ${\bf P}$ is locally bounded and smooth in $\mathbb{R}^2\setminus\Sigma$, the function $u$ is a classical solution of (\ref{eq 2.2}) in $\mathbb{R}^2\setminus\Sigma$. By standard regularity theory it belongs to $C^\infty\left(\mathbb{R}^2\setminus\Sigma\right)$.
Then $w$ is a smooth locally bounded function in $\mathbb{R}^2\setminus\Sigma$ satisfying (\ref{eq 3.2}), an equation that we rewrite under the form
\begin{equation}\label{eq 3.2'}
\displaystyle - \Delta w +h(\cdot,w)=f_1-f_2\quad{\rm in}\;\; \CD'(\R^2),
\end{equation}
where the function 
$h(x,z)$ is defined in $\BBR^2\times\BBR$ by 
 $$h(x,z)=\arraycolsep=1pt\left\{
\begin{array}{lll}
 \displaystyle {V}(x)\frac{  e^{z}}{(e^{\nu_1-\nu_2}+  e^{z})^{1+a}}\qquad  &{\rm for }\ \ x\in \mathbb{R}^2\setminus\Sigma,
\\[2.5mm]\phantom{ }
 \displaystyle 0\qquad    &{\rm for }\ \ x\in\Sigma_2,
\\[2.5mm]\phantom{ }
 \displaystyle  \gs^{-2an_j}\prod_{i\neq k}|p_j-p_i|^{-2an_i}e^{-az} \quad  &{\rm for }\ \   x=p_j\in \Sigma_1.
\end{array}
\right.$$
The function $h$ is nonnegative and smooth in $\R^2\setminus\Sigma$ and continuous in $\BBR^2\times\BBR$.
Since $w$ is smooth in $\R^2\setminus\Sigma$, so is $h(\cdot,w)$. Next we set, with $Z=e^z\geq 0$
\begin{equation}\label{eq 3.2+1}\phi(Z)=\frac{Ze^{a(\nu_1-\nu_2)}}{\left(e^{a(\nu_1-\nu_2)}+Z\right)^{1+a}}\Longrightarrow \phi'(Z)=\frac{e^{a(\nu_1-\nu_2)}\left(e^{a(\nu_1-\nu_2)}-aZ\right)}{\left(e^{a(\nu_1-\nu_2)}+Z\right)^{2+a}}.
\end{equation}
Then 
\begin{equation}\label{eq 3.2+2}\phi'(Z_0)=0\quad\text{with }
\; Z_0=\frac{e^{a(\nu_1-\nu_2)}}{a}\Longrightarrow \phi(Z_0)=\frac{a^a}{(a+1)^{1+a}}e^{(a-a^2)(\nu_1-\nu_2)}=\max\{\phi(Z):Z>0\}.
\end{equation}
Hence 
\begin{equation}\label{eq 3.2+3}
0\leq h(x,w)\leq {\bf P}(x)\frac{a^a}{(a+1)^{1+a}}e^{(a-a^2)(\nu_1-\nu_2)}. 
\end{equation}
Note that ${\bf P}$ is locally bounded in $ \mathbb{R}^2\setminus\Sigma_1$, then it follows by standard regularity arguments, (see e.g. \cite{G}) that 
$w$ belongs to $W^{2,t}_{loc}( \mathbb{R}^2\setminus\Sigma_1)$ for any $1<t<\infty$. Hence $h(\cdot,w)\in C^{1,\gth}(\mathbb{R}^2\setminus\Sigma_1)$ for any $\gth\in (0,1)$, and finally $w\in C^{3,\gth}(\mathbb{R}^2\setminus\Sigma_1)$ is a strong solution in $\mathbb{R}^2\setminus\Sigma_1$. In a neighborhood of $\Sigma_1$ we write $h$ under the form
$$h(x,z)= {\bf P}(x)e^{a\nu_1}\frac{  e^{z-a\nu_2}}{(e^{\nu_1-\nu_2}+  e^{z})^{1+a}}.
$$
Since $h$ is nonnegative, then $w$ satisfies the inequality 
$$-\Gd w\leq f_1-f_2\quad\text{in }\;\CD'(\R^2),
$$
and as $f_1-f_2$ is bounded with compact support, it follows that $w$  is locally bounded from above in $\R^2$.
Furthermore, there exist an open set $\CO$ such that $\Gs_1\subset\CO$ and $\overline\CO\cap \Gs_2=\emptyset$ and a function $\zeta_1\in C(\overline \CO)$  such that 
$$h(x,z)=\zeta_1 \frac{e^z}{(e^{\nu_1-\nu_2}+  e^{z})^{1+a}}\quad\text{for all }\; (x,z)\in \overline \CO\times \R.$$
For a given $p_j\in\Gs_1$, we set $r_j=\sup\{w(x):x\in\overline B_\gs(p_j)\}$ and $v_j=r_j-w$. Then $v_j\geq 0$ in $\overline B_\gs(p_j)$ and 
$$-\Gd v_j=f_2-f_1+\zeta_1\frac{e^{r_j-v_j}}{(e^{\nu_1-\nu_2}+  e^{r_j-v_j})^{1+a}}=f_2-f_1+\zeta_1e^{-ar_j}\frac{e^{-v_j}}{(e^{\nu_1-\nu_2-r_j}+  e^{-v_j})^{1+a}}.$$
Since $f_1,f_2$ are smooth, hence  $\zeta_1e^{-ar_j}\frac{e^{-v_j}}{(e^{\nu_1-\nu_2-r_j}+  e^{-v_j})^{1+a}}\in L^1(B_\gs(p_j))$ by \cite{BrLi}. 
If $0<\gs'\leq\gs$, we denote by $\phi^{B_{\gs'}}_j$ the harmonic lifting of $v_j\lfloor_{\partial B_{\gs'}}$ in $B_{\gs'}(p_j)$ and put $\tilde v_{\gs'}=v_j-\phi^{B_{\gs'}}_j$. Then for $\gs'\leq\gs$,
$$\left\{\BA{lll} -\Gd \tilde v_{\gs'}=f_2-f_1+\zeta_1e^{-ar_j}\frac{e^{-v_j}}{(e^{\nu_1-\nu_2-r_j}+  e^{-v_j})^{1+a}}:=F_j\quad&\text{in }\;\, B_{\gs'}(p_j),\\[2mm]
\phantom{-\Gd }\tilde v_{\gs'}=0\quad&\text{on }\;\partial B_{\gs'}(p_j).
\EA\right.$$
Let $M^2(B_{\gs'}(p_j))$ denote the Marcikiewicz space also known as the Lorentz space $L^{2,\infty}(B_{\gs'}(p_j))$. Then there holds
\begin{equation}\label{eq 3.2+4}
\norm{\nabla \tilde v_{\gs'}}_{M^2(B_{\gs'}(p_j))}\leq c_0\norm{F_j}_{L^1(B_{\gs'}(p_j))}
\end{equation}
and the constant $c_0$ is independent of $\gs'$. We recall below John-Nirenberg's theorem \cite[Theorem 7.21]{G}: 
{\it Let $u\in W^{1,1}(G)$ where $G\subset\Gw$ is convex and suppose there is a constant $K$ such that 
\begin{equation}\label{eq 3.2+5}
\int_{G\cap B_r}|\nabla u|dx\leq Kr\quad\text{for any ball }\, B_r(0).
\end{equation}
Then there exist positive constants $\gm_0$ and $c_1$ such that 
\begin{equation}\label{eq 3.2+6}
\int_{G}\exp\left(\frac{\gm}{K}|u- u^{_G}|\right) dx\leq c_1(diam (G))^2,
\end{equation}
where $\gm=\gm_0|G|(diam (G))^{-2}$ and $ u^{_G}$ is the average of $u$ on $G$.
}\smallskip

\noindent From (\ref{eq 3.2+4}) with $G=B_{\gs'}(p_j)$,
\begin{equation}\label{eq 3.2+7}
\int_{B\cap B_r}|\nabla \tilde v_{\gs'}|dx\leq c_3r\norm{F_j}_{L^1(B_{\gs'}(p_j))}:=K(\gs')r,
\end{equation}
 and since $|G|(diam (G))^{-2}=\gp$, we obtain
 \begin{equation}\label{eq 3.2+8}
\int_{B_{\gs'}(p_j)}\exp\left(\frac{\gp\gm_0}{K(\gs')}|\tilde v_{\gs'}-\tilde v_{\gs'}^{_{B_{\gs'}}}|\right) dx\leq c_1\gs'^2.
\end{equation}
 Hence, for any $\gk>0$ there exists $\gs'\in (0,\gs]$ such that 
  \begin{equation}\label{eq 3.2+8+}
\int_{B_{\gs'}(p_j)}\exp\left(\gk|\tilde v_{\gs'}-\tilde v_{\gs'}^{_{B_{\gs'}}}|\right) dx\leq c_1\gs'^2\Longrightarrow 
\int_{B_{\gs'}(p_j)}\exp(\gk\tilde v_{\gs'})dx)\leq c_1\gs'^2\exp(\gk\tilde v_{\gs'}^{_{B_{\gs'}}}).
\end{equation}
Now we observe that there holds in $B_{\gs'}(p_j)$,
$$|F_j|\leq |f_2-f_1|+\zeta_1e^{-ar_j}e^{av_j}\leq |f_2-f_1|+\zeta_1e^{a\sup|v_j\lfloor_{\partial B_{\gs'}}|}e^{a\tilde v_{\gs'}}.
$$
For $\gk>a$, 
$$\int_{B_{\gs'}(p_j)}|F_j(x)|^{\frac{\gk}{a}}dx \leq 2^{\frac{\gk}{a}-1}\int_{B_{\gs'}(p_j)}\left(|f_2-f_1|^{\frac{\gk}{a}}
+\left(\zeta_1e^{a\sup|v_j\lfloor_{\partial B_{\gs'}}|}\right)^{\frac{\gk}{a}}e^{\gk\tilde v_{\gs'}}\right) dx.
$$
By (\ref{eq 3.2+8}) the right-hand side of the above inequality is bounded, hence $F_j\in L^{\frac{\gk}{a}}(B_{\gs'}(p_j))$. Since $\tilde v_{\gs'}$ vanishes on $\partial B_{\gs'}(p_j)$, it follows by $L^p$ regularity theory that $\tilde v_{\gs'}\in W^{2,\frac{\gk}{a}}(B_{\gs'}(p_j))\cap W^{1,\frac{\gk}{a}}_0(B_{\gs'}(p_j))$. By Sobolev embedding theorem, $\tilde v_{\gs'}\in L^\infty(B_{\gs'}(p_j))$. Hence $F_j\in L^\infty(B_{\gs'}(p_j))$ and again 
$\tilde v_{\gs'}\in W^{2,q}(B_{\gs'}(p_j))$ for any $q\in [1,\infty)$ and thus $\tilde v_{\gs'}\in C^{1,\gth}(\overline B_{\gs'}(p_j))$ for any $\gth\in (0,1)$. Therefore $v_j$ remains bounded in $C^{1,\gth}(\overline B_{\gs''}(p_j))$ for any $\gs''<\gs'$. In a neighborhood of $p_j$, 
$x\mapsto |x-p_j|^{-2an_j}e^{2a\gz(\frac{ |x-p_j|}{\gs})}$ is H\"older continuous (of order $2an_j$ if $2an_j<1$), and so is  
$x\mapsto {\bf P}(x)e^{a(\gn_1-\gn_2)(x)}$. For the same reason, $x\mapsto \frac{e^{w(x)}}{\left(e^{(\gn_1-\gn_2)(x)+e^{w(x)}}\right)^{1+a}}$ is H\"older continuous (with the same exponent) near $p_j$. Finally we infer that there exists $\gth\in (0,1)$ such that $v_j\in C^{2,\gth}(\overline B_{\gs''}(p_j))$, which implies that $w\in C^{2,\gth}(\BBR^2)$ is a strong solution of (\ref{eq 3.2}) in $\BBR^2$.
\hfill$\Box$\medskip

\noindent{\it Remark.} Since $\gv_1$ and $\gn_2$ have compact support, we note that a weak solution $ u_\beta$ with the asymptotic behavior
 $\beta\ln |x|+O(1)$ at infinity can be decomposed 
 \begin{equation}\label{w beta}
  u_\beta=w_\beta-\nu_1+\nu_2,
 \end{equation}
where $w_\beta$ is a classical solution of (\ref{eq 3.2}) with the same asymptotic behavior
 $\beta\ln |x|+O(1)$ at infinity. In fact, we shall continue to take out the singular source of the solution $w_\beta$ at infinity in our derivation
 of non-topological solutions of (\ref{eq 1.1}).

\subsection{Basic estimates }
The following  estimates play an important role in our construction of solutions to (\ref{eq 1.1}).

\begin{lemma}\label{lm 2.1}
Let $\Gamma$ be the fundamental solution of $-\Delta$ in $\R^2$, $F\in L^p_{loc}(\R^2)$, $p>1$, with the support in $\overline{B_R(0)}$ for some $R>0$ such that
\begin{equation}\label{2.2}
\int_{\R^2} F(x)dx=0.
\end{equation}
Then there holds
 \begin{equation}\label{e 2.1+1}
|\Gamma\ast F(x)|\le \frac{R}{|x|} \norm{F}_{L^1(\R^2)}\quad{\rm for}\  |x|>4R,
\end{equation}
and for some $c_2>0$ depending on $p$ and $R$,
\begin{equation}\label{e 2.3}
\norm{\Gamma\ast F}_{L^\infty(\R^2)}\le   c_2 \norm{F}_{L^p(\R^2)}.
\end{equation}

\end{lemma}
{\bf Proof. } As supp$(F)\subset \overline{B_R(0)}$ and $F\in L^p_{loc}(\R^2)$,  $F\in L^1(\R^2)\cap L^p(\R^2)$ by H\"older's inequality. Since (\ref{2.2}) holds, we have for $ |x|>4R$,
\begin{eqnarray*}
|\Gamma\ast F(x)|&=& \frac{1}{2\pi}\Big|\int_{B_R(0)}  \ln|x-y|\, F(y)dy - \int_{B_R(0)}  \ln|x| \,F(y)dy\Big| \\
&= & \frac{|x|^2}{2\pi}\Big|\int_{B_\frac{R}{|x|}(0)}   \ln|e_x-z|\,  F(|x|z)dz\Big|
\\& \le &  \frac{|x|^2}{\pi} \int_{B_\frac{R}{|x|}(0)} |z| |F(|x|z)|dz
 \\&\le & \frac {2R}{\pi |x|}   \int_{B_ {R} (0)}  |F(y)|dy
\\& < & \frac{R}{|x|} \norm{F}_{L^1(\R^2)},
\end{eqnarray*}
where $e_x=\frac{x}{|x|}$. Besides (\ref{2.2}) we have used  the fact that
$$|\ln|e_x-z||\le 2|z|\le 2 \frac{R}{|x|}\quad {\rm for\ any}\  \, z\in B_{ {R}/{|x|}}(0)\subset B_{1/4}(0).$$
Therefore, (\ref{e 2.1+1}) is proved. On the other hand, for $|x|\le 4R$, we have that
\begin{eqnarray*}
|\Gamma\ast F(x)|&=&\frac1{2\pi}\Big|\int_{B_R(0)} F(y)\ln|x-y| dy\Big|
\\&\le&  \Big(\int_{B_R(0)} |F(y)|^p dx\Big)^{\frac1p} \Big(\int_{B_R(0)} |\ln |x-y||^{p'} dx\Big)^{\frac1{p'}}
\\&\le & c_p\norm{F}_{L^p(B_R(0))},
\end{eqnarray*}
where $p'=\frac{p}{p-1}$ and $\displaystyle c_p=\max_{|x|\le 4R}(\int_{B_R(0)} |\ln |x-y||^{p'} dx)^{\frac1{p'}} $.
Thus, (\ref{e 2.3}) follows and   the proof is complete.
 \hfill$\Box$\medskip

For functions with non-compact supports, we have the following estimates.

\begin{lemma}\label{lm 1.2}
Let $F\in L^p_{loc}(\R^2)$ with $p>1$ satisfy that
\begin{equation}\label{2.2-9-5}
  \int_{\R^2} F(x)\,dx =0
\end{equation}
and
\begin{equation}\label{2.4-9-5}
   |F(x)|\le c_3|x|^{-\tau}\quad {\rm for}\  |x|\geq r
\end{equation}
for some $\tau>2$, $c_3>0$ and $r>0$. Then  for some $c_4>0$
$$\norm{\Gamma\ast F}_{L^\infty(\R^2)}\le c_4,$$
and there exist  $ c_5>0$ and $r_0>r$ such that for $|x|\geq r_0$
\begin{equation}\label{2.3-9-5}
|x||\nabla\Gamma\ast F(x)|+|\Gamma\ast F(x)|\le   \frac{c_5}{(\tau-2)^2} |x|^{-\frac{\tau-2}{\tau-1}}.
\end{equation}

\end{lemma}
{\bf Proof.} If $F\in L^p_{loc}(\BBR^2)$ satisfies (\ref{2.4-9-5}), then $F\in L^1(\BBR^2)$. Let $\eta_r:\R^2\to [0,1]$ be a smooth and radially symmetric function such that
$$\eta_r=1\ \, {\rm in}\  B_r(0),\ \quad \eta_r=0\ \, {\rm in}\ \, B_{r+1}(0),$$
and denote
$$F_1=F\eta_r-\left(\int_{\BBR^2}F\eta_r dx\right)\frac{\eta_r}{\norm{\eta_r}_{L^1(\R^2)}},\ \quad F_2=F-F_1.$$
By (\ref{2.2-9-5}), we have that
$$\int_{\R^2} F_1\, dx=\int_{\R^2} F_2\, dx=0.$$
Since $F_1\in L^p_{loc}(\R^2)$  has compact support, it follows by Lemma \ref{lm 2.1},
$\Gg\ast F_1$ is bounded  and satisfies (\ref{e 2.1+1}).
Concerning $F_2$, we have
$$F_2=F-F_1=\frac{\int_{\BBR^2}F\eta_r dx}{\norm{\eta_r}_{L^1(\R^2)}}\quad\text{on }\, B_r,
$$
and $F_2$ satisfies (\ref{2.4-9-5}) on $B_r^c$, with may be another constant. It is locally bounded hence
$\Gamma\ast F_2$ is also locally bounded  in $\R^2$. If  (\ref{2.4-9-5}) holds true, we have to prove that $F_2$ verifies (\ref{2.3-9-5}).
\smallskip

\noindent Since $\int_{\R^2} F_2\, dx=0$, then for all $|x| > 4 r$ and $R\in (r,\frac{|x|}4)$ which will be chosen latter on,
\begin{eqnarray*}
2\pi(\Gamma\ast F_2)(x)&=&   |x|^2 \int_{\R^2}   \ln|e_x-z| F_2(|x|z)dz +  |x|^2\ln |x|\int_{\R^2}   F_2(|x|z)dz
\\&=&   |x|^2 \int_{B_{R/|x|}(0)}  \ln|e_x-z|  F_2(|x|z)dz +|x|^2 \int_{ B_{1/2}(e_x)}  \ln|e_x-z|  F(|x|z)dz
\\&&+|x|^2 \int_{\R^2\setminus (B_{R/|x|}(0)\cup B_{1/2}(e_x))}  \ln|e_x-z| \, F_2(|x|z)dz
\\&=: & I_1(x)+I_2(x)+I_3(x),
\end{eqnarray*}
using the fact that $B_{R/|x|}(0)\cap B_{1/2}(e_x))=\emptyset$. By a direct computation, we have that
$$\begin{array}{lll}\displaystyle
|I_1(x)| \le    |x|^2 \int_{B_{R/|x|}(0)} |z| |F_2(|x|z)|dz
\\[4mm]\displaystyle
\phantom{|I_1(x)| }=   2\frac{R}{|x|} \int_{B_{R}(0)} | F_2(y)|dy 
\\[4mm]\displaystyle
\phantom{|I_1(x)| } \leq 2\frac{R}{|x|} \norm{F_2}_{L^1(\R^2)}.
\end{array}
$$
For $z\in B_{1/2}(e_x)$, there holds $|x||z|\geq \frac12|x|>2r$, then $|F(|x|z)|\le c_3|x|^{-\tau}|z|^{-\tau}$ and
$$\begin{array}{lll}\displaystyle
|I_2(x)| \leq   c_3 |x|^{2-\tau} \int_{B_{1/2}(e_x)} (-\ln|e_x-z|) |z|^{-\tau} dz
\\[4mm]\displaystyle \phantom{|I_2(x)| }\leq   2^\tau c_3 |x|^{2-\tau} \int_{B_{1/2}(e_x)} (-\ln|e_x-z|)   dz
\\[4mm]\displaystyle\phantom{|I_2(x)| }\leq c_6 R^{2-\tau},
\end{array}
$$
where $c_6 =  2^{2(N-M)}c_3 \big(\int_{B_{1/2}(0)} (-\ln|z|)   dz\big)$ can be chosen  independently of $\tau $ in $ (2, 2(N-M))$.\\
Next, if $z\in \R^2\setminus (B_{R/|x|}(0)\cup B_{1/2}(e_x))$, then $\big|\ln |e_x-z|\big|\le \ln (1+|z|)$ and $|F(|x|z)|\le c_7|x|^{-\tau}|z|^{-\tau}$, since $|z| \geq \frac{R}{|x|} > \frac{r}{|x|}$. By the integration by parts we get
\begin{eqnarray*}
|I_3(x)|\!\!\!\! &\le &\!\!\!\!  c_8  |x|^{2-\tau} \int_{\R^2\setminus  B_{R/|x|}(0)  } \ln(1+|z|)\, |z|^{-\tau}dz
\\\!\!\!\!&\le &\!\!\!\!  \frac{2\pi c_8}{\tau-2} R^{2-\tau} \ln \left(1+\frac{R}{|x|}\right) + \frac{2\pi c_5}{(\tau-2)^2} R^{2-\tau}
\\\!\!\!\!&\le &\!\!\!\!   \frac{2\pi c_8}{(\tau-2)^2}\Big((\tau-2)\ln 2+1 \Big)R^{2-\tau}.
\end{eqnarray*}
Thus, taking $R=|x|^{\frac1{\tau-1}}$ and $|x|$ sufficiently large (certainly $R\in (r,\frac{|x|}4)$ is satisfied), we have
\begin{eqnarray*}
|\Gamma\ast F_2(x)| &\le &   \frac{R}{\pi |x|} \norm{F}_{L^1(\R^2)} +\frac{c_8}{2\pi}  R^{2-\tau} +\frac{ c_8}{(\tau-2)^2}\Big(2(N-M-1)\ln (e+1)+1\Big) R^{2-\tau}
\\&\le &   \frac{c_9}{(\tau-2)^2} |x|^{-\frac{\tau-2}{\tau-1}},
\end{eqnarray*}
where $c_9>0$ can be chosen  independently of $\tau$. In order to prove the gradient estimate, we denote by $(r,\gth)$ the polar coordinates in 
$\BBR^2$, set $t=\ln r$ and 
$$\gw(t,\gth)=\tilde \gw(r,\gth)=r^{-\frac{\gt-2}{\gt-1}}\Gamma\ast F(r,\gth)\quad\text{and }\; \gf(t,\gth)=r^{\gt}F(r,\gth).
$$
Then $\gw$ and $\gf$ are bounded on $[\ln r_1,\infty)\times S^1$ where there holds
$$\CL\gw:=\frac{\partial^2\gw}{\partial t^2}-2\frac{\gt-2}{\gt-1}\frac{\partial\gw}{\partial t}+\left(\frac{\gt-2}{\gt-1}\right)^2\gw+\frac{\partial^2\gw}{\partial \gth^2}
=\gt^{-\frac{\left(\gt-2\right)^2}{\gt-1}}\gf. 
$$
Since the operator $\CL$ is uniformly elliptic on $[\ln r_1,\infty)\times S^1$, for any $T>\ln r_1+2$ there holds by standard elliptic equations regularity estimates \cite{G},
$$\BA{lll}\displaystyle\sup_{[T-1,T+1]\times S^1}\left(\left|\frac{\partial\gw}{\partial \gth}\right|+\left|\frac{\partial\gw}{\partial t}\right|\right)
\leq c_{10}\displaystyle\sup_{[T-2,T+2]\times S^1}\left(\left|\gw\right|+\left|\gt^{-\left(\frac{\gt-2}{\gt-1}\right)^2}\gf\right|\right)
\leq c_{11},
\EA$$
and $c_{11}$ does not depend on $T$. As $|x||\nabla \tilde\gw(x)|=\left(\left|\frac{\partial\gw}{\partial \gth}\right|^2+\left|\frac{\partial\gw}{\partial t}\right|^2\right)^\frac{1}{2}$, this implies the claim. 
  \hfill$\Box$\medskip
  
  The decay estimate on the gradient at infinity does not use the  identity (\ref{2.2-9-5}). It is actually more general.
  
\begin{corollary}\label{gradecay} Let  $F\in L^p_{loc}(\R^2)$ with $p>1$ satisfy (\ref{2.4-9-5}) with $\gt>2$ and $w$ be a solution of 
\begin{equation}\label{Z1} 
-\Gd w=F\quad\text{in }\;\BBR^2.
\end{equation}
(i) If $\displaystyle\lim_{|x|\to\infty}|x|^{\frac{\gt-2}{\gt-1}}|w(x)|<\infty $, then
\begin{equation}\label{Z2} 
\abs{\nabla w(x)}\leq c _{12}|x|^{-\frac{2\gt-3}{\gt-1}}\quad\text{for }\; |x|\geq r_0.
\end{equation}
(ii) If there exists a constant $C$ such that $w(x)=C+O(|x|^{-\frac{\gt-2}{\gt-1}})$ when $|x|\to\infty$, then estimate (\ref{Z2}) holds.
\end{corollary}
{\bf Proof.} The assertion (i) is clear since the starting point of the gradient estimate in the previous lemma is 
$$|w(x)|\leq c_{13}|x|^{-\frac{\gt-2}{\gt-1}}\quad\text{for $|x|$ large enough}.
$$
For assertion (ii), we set $w(x)=C+\tilde w(x)$ where $|\tilde w(x)|=O(|x|^{-\frac{\gt-2}{\gt-1}})$. Then $-\Gd\tilde w=-\Gd w$ and $\nabla w=\nabla \tilde w$. We conclude by (i).   \hfill$\Box$\medskip

When $\gt=2$, Lemma \ref{lm 1.2} is no longer valid, however the following limit case is available. 
\begin{lemma}\label{lm 1.3}
Let $F\in L^p_{loc}(\R^2)$ with $p>1$ satisfy (\ref{2.2-9-5})
and
\begin{equation}\label{3.4-9-5}
   |F(x)|\le c_{14}|x|^{-2}(\ln |x|)^{-\gn}\quad {\rm for}\quad |x|\geq r,
\end{equation}
for some $\gn>2$, $c_{14}>0$ and $r>0$. Then
$$\norm{\Gamma\ast F}_{L^\infty(\R^2)}\le c_{15},$$
and 
\begin{equation}\label{3.3-9-5}
|x||\nabla\Gamma\ast F(x)|+|\Gamma\ast F(x)|\le   c_{16}(\ln |x|)^{-\gn}\quad\text{for \ $|x|\geq r_1$},
\end{equation}
where $c_{15},\, c_{16}>0$ and $r_1> r$ is large enough.

\end{lemma}
{\bf Proof}. The assumption (\ref{3.4-9-5}) jointly with $F\in L^p_{loc}(\R^2)$ implies $F\in L^1(\R^2)$. We write $F=F_1+F_2$  in the same way as in Lemma \ref{lm 1.2}. Clearly $\Gg\ast F_1$ is uniformly bounded and satisfies (\ref{e 2.1+1}) and
 (\ref{e 2.3}). Then for $|x|>4r>4$ and $R\in (r,\frac{|x|}{4}]$, 
 $$\BA {lll}2\gp(\Gg\ast F_2)(x)=I_1(x)+I_2(x)+I_3(x),
 \EA$$
 where $I_1$, $I_2$ and $I_3$ are defined in the proof of Lemma \ref{lm 1.2} and where $|I_1(x)|\leq 2\frac{R}{|x|}\norm{F_2}_{L^1(\BBR^2)}$. 
 
 When $z\in B_{\frac 12}(e_x)$ we have $|x||z|\geq \frac 12|x|>2r$, hence 
\begin{equation}\label{3.3-9-6}|F_2(|x|z)|\leq c_{14}|x|^{-2}|z|^{-2}\left(\ln |x|+\ln|z|\right)^{-\gn}\leq c_{14}|x|^{-2}|z|^{-2}\left|\ln |x|-\ln2\right|^{-\gn},\end{equation}
 and
 $$\BA {lll}\displaystyle
 I_2(x)\leq -c_{14}\left|\ln |x|-\ln2\right|^{-\gn}\int_{B_{\frac 12}(e_x)}\ln |z-e_x||z|^{-2} dz\leq c_{17}\left(\ln |x|\right)^{-\gn}.
 \EA$$
Finally, if $z\in \R^2\setminus (B_{R/|x|}(0)\cup B_{1/2}(e_x))$, then $|\ln |e_x-z||\le \ln (1+|z|)$ and 
$$\BA{lll}|F_2(|x|z)|\le c_{17}(1+|x||z|)^{-2}\left(\ln (1+|x||z|)\right)^{-\gn}\le c_{17}|x|^{-2}|z|^{-2}\left(\ln (1+|x||z|)\right)^{-\gn}.
\EA$$ 
Since $|z| \geq \frac{R}{|x|} > \frac{r}{|x|}$, we have
$$\BA{lll}\displaystyle
I_3(x)\leq c_{18}\int_{\frac {R}{|x|}}^{\infty}\ln (1+t)(\ln (1+t|x|)^{-\gn}\frac{dt}{t}\leq c_{18}\int_{R}^{\infty}(\ln (1+s))^{1-\gn}\frac{ds}{s}.
\EA$$ 
Since $R>r>1$, 
$$\BA{lll}\displaystyle
I_3(x)\leq c_{19}\int_{R}^{\infty}(\ln (1+s))^{1-\gn}\frac{ds}{1+s}=\frac{c_{19}}{\gn(\ln (1+R))^{\gn}}.
\EA$$ 
If we choose $R=\frac{|x|}{4}$, we obtain that  $(\ln (1+|x|))^{\gn}(\Gg\ast F)(x)$ remains uniformly bounded on $\BBR^2$. \smallskip

 Next we prove the gradient estimate. Set $t=\ln r$, $\gw(t,\gth)=(\Gg\ast F)t,\gth)$ and $\gf(t,\gth)=F(r,\gth)$, then 
$$\CL\gw:=\frac{\partial^2\gw}{\partial t^2} +\frac{\partial^2\gw}{\partial \gth^2}
=\gf
$$
and $|\gw(t,\gth)|\leq c_{20}t^{-\gn}$ and $|\gf(t,\gth)|\leq c_{20}t^{-\gn}$ for $t\geq t_1$. Since the operator $\CL=\frac{\partial^2 }{\partial t^2}+\frac{\partial^2 }{\partial \gth^2}$ is uniformly elliptic, then we have for $T\geq \max\{4,t_1\}$,
$$\BA{lll}\displaystyle \sup_{[T-1,T+1]\times S^1}\left(\left|\frac{\partial\gw}{\partial \gth}\right|+\left|\frac{\partial\gw}{\partial t}\right|\right)
\leq c_{10}\displaystyle\sup_{[T-2,T+2]\times S^1}\left(\left|\gw\right|+\left|\gf\right|\right)
\leq 2c_{10}c_{20}(T-2)^{-\gn}\leq  c_{21}T^{-\gn}.
\EA$$
Returning to the variable $x$, we infer (\ref{3.3-9-5}).    \hfill$\Box$\medskip

Similarly as in Corollary \ref{gradecay}, the following extension of (\ref{3.3-9-5}) holds.

\begin{corollary}\label{gradecay*} Let  $F\in L^p_{loc}(\R^2)$ with $p>1$ satisfy (\ref{3.4-9-5}) with $\gn>2$ and $w$ satisfiesf (\ref{Z1}).\smallskip

\noindent (i) If $\displaystyle\lim_{|x|\to\infty}(\ln |x|)^\gn |w(x)|<\infty$, then there exist $c_{22}>0$ and $r_2>1$ such that 
\begin{equation}\label{Z2'} 
\abs{\nabla w(x)}\leq   c_{22}|x|^{-1}(\ln |x|)^{-\gn}\quad\text{for }\; |x|\geq r_2.
\end{equation}
(ii) If there exists a constant $C$ such that $w(x)=C+O((\ln |x|)^{-\gn})$ when $|x|\to+\infty$, then estimate (\ref{Z2'}) holds.
\end{corollary}

\subsection{Related problems with increasing nonlinearity  }

In order to remove the  condition $\beta\ln|x|+O(1)$ as $|x|\to\infty$ satisfied by the solutions of (\ref{eq 1.1}), we introduce two functions $\lambda$ and  $\Lambda$, which are   positive smooth functions such that
\begin{equation}\label{lambda}
\lambda(x)=|x|,\quad \Lambda(x)= \ln |x|\quad {\rm for}\ |x|\geq e^e.
\end{equation}
Since $\Delta \Lambda=0$ in $B_{e^e}^c(0)$, 
\begin{eqnarray}
\Delta\ln\Lambda =  \frac{\Delta \Lambda}{\Lambda}-\frac{ |\nabla \Lambda|^2}{\Lambda^2} =   -\frac{1}{|x|^2(\ln |x|)^2}\quad {\rm in}\ \, B_{e^e}^c(0), \label{e 2.5}
\end{eqnarray}
and
\begin{equation}\label{Vl}
\frac1{2\pi}\int_{\R^2}(\Delta \ln \lambda) dx =1,
\end{equation}
it implies 
\begin{eqnarray*}
\frac1{2\pi}\int_{\R^2}(\Delta \ln \Lambda) dx  &=& \lim_{r\to+\infty}\frac1{2\pi}\int_{\partial B_r(0)} \frac{\nabla\Lambda(x) }{\Lambda(x)}\cdot \frac{x}{|x|} d\omega(x) \\
     &=&    \lim_{r\to+\infty}\frac{1}{r\ln r}
   = 0.
\end{eqnarray*}
In what follows we classify the  solutions of the following equations
\begin{equation}\label{eq 2.1}
-\Delta u +  W F_i(\lambda^\beta e^u) =g_\beta \quad{\rm in}\;\;  \mathbb{R}^2,
\end{equation}
where $i=1,2$, $F_1(s)=s$,  $F_2(\cdot,s)=\frac{s}{e^{(\nu_1-\nu_2)(\cdot)}+s}$,
  \begin{equation}\label{2.5}
 g_\beta =f_1-f_2+\beta\Delta \ln \lambda,
 \end{equation}
 and where  $W$  satisfies the following assumption:\smallskip

\noindent $(\mathcal{W}_0)$ {\it
The function $W$  is positive and locally H\"{o}lder continuous in $\R^2\setminus \Sigma_1$ and
$$W(x)\le c_{23}|x-p_j|^{-\tau_{p_j}}\quad\text{in}\quad B_{\sigma}(p_j)\ \quad{\rm and}\quad\ \limsup_{|x|\to+\infty} W(x)|x|^{\gamma_\infty}<+\infty, $$
\qquad\ where $c_{23}>0$, $  \tau_{p_j}\in[0,2)$ and $\gamma_\infty>0$.}
\smallskip

It is important to note that from (\ref{Vf}), (\ref{Vl}) and (\ref{2.5}), there holds
\begin{equation}\label{gbeta}\int_{\R^2} g_\beta\, dx=2\pi[2(N-M)+ \beta]. \end{equation}
\begin{teo}\label{teo 2.1}
Assume that  $F_1(s)=s$,  $F_2(x,s)=\frac{s}{e^{(\nu_1-\nu_2)(x)}+s}$  and  $g_\beta$ is a  H\"{o}lder continuous function   with compact support in $B_{r_0}(0)$ satisfying the relation (\ref{gbeta}) 
for some nonnegative integers $N$ and $M$. Let $W$ verify $(\mathcal{W}_0)$ with
 \begin{equation}\label{2La}
\gg_\infty>\gb+2.
 \end{equation}

\noindent $(i)$  
Then problem (\ref{eq 2.1}) with $i=1$ has a unique bounded solution $v$ verifying, for some $C_\beta\in\R$,
\begin{equation}\label{2.1 e}\BA{lll}
v(x)=C_\beta+O(|x|^{-\frac{\gamma_\infty-\beta-2 }{\gamma_\infty-\beta-1 }})\,\text{ and }\;
  |\nabla v(x)|=O(|x|^{-1-\frac{\gamma_\infty-\beta-2 }{\gamma_\infty-\beta-1 }})\quad\text{as }\, |x|\to+\infty.
\EA\end{equation}
$(ii)$  Assume additionally that $2(N-M)+ \beta>0$ and
\begin{equation}\label{2.1-1}
2\pi[2(N-M)+ \beta]< \int_{\R^2}W dx\leq+\infty.
\end{equation}
Then  problem (\ref{eq 2.1}) with $i=2$ has a unique bounded solution verifying 
\begin{equation}\label{2.1 e2}\BA{lll}
v(x)=C_\beta+O(|x|^{-\frac{\gamma_\infty+\beta_--2 }{\gamma_\infty+\beta_--1 }})\,\text{ and }\;
  |\nabla v(x)|=O(|x|^{-1-\frac{\gamma_\infty+\beta_--2 }{\gamma_\infty+\beta_--1 }})\quad\text{as }\, |x|\to+\infty.
\EA\end{equation} 
\end{teo}
{\bf Proof}.   
 {\it Step 1}. Since $0\leq\gt_{p_j}<2$, $W\in L^1_{loc}(\BBR^N)$. For $t\in\R$, we set
$$h_{i,t}(x)=   W(x)F_i(x,\lambda^\beta(x) e^t)\qquad  \forall\,   x\in \R^2\setminus\Sigma_1.$$
Notice that $h_{2,t}$ is defined on $\Gs_2$ by 
$$h_{2,t}(q_j)=0=\lim_{x\to q_j}h_{2,t}(x)\quad\text{for all }\, q_j\in\Gs_2.
$$
The function $h_{i,t}$ is H\"{o}lder continuous in $\R^2\setminus \Sigma_1$,
  $t\mapsto h_{i,t}$ is  increasing in $\R^2\setminus\Sigma_1$, and there holds
 $$h_{1,t}\to \infty\quad{\rm locally\ in}\ \R^2\setminus \Sigma_1 \quad {\rm as}\ t\to+\infty,$$
$$h_{2,t}\to W\quad{\rm locally\ in}\ \R^2\setminus \Sigma_1 \quad {\rm as}\ t\to+\infty.$$
Furthermore,
 $$h_{i,t}\to 0\quad{\rm locally\ in}\ \R^2\setminus \Sigma_1  \quad {\rm as}\ t\to-\infty,\quad i=1,2.$$
Using assumption $(\mathcal{W}_0)$, we obtain that
\begin{equation}\label{2.2a}h_{1,t}(x)\leq c_{24}e^t|x|^{-\gg_\infty+\gb}\quad\text{for }\,|x|\geq r_3
\end{equation}
for some $r_3>0$. Since $-\gg_\infty+\gb<-2$, we have that 
$$\lim_{t\to-\infty}\int_{\R^2}h_{1,t}(x)dx=0.
$$
Concerning $h_2$, we have $(\gn_1-\gn_2)(x)(x)=0$ if dist$ (x,\Gs)\geq\gs$, and there holds 
$$F_2(x,\lambda(x)^\beta e^t)=\frac{\lambda(x)^\beta e^t}{1+ \lambda(x)^\beta e^t}=\frac{|x|^\beta e^t}{1+ |x|^\beta e^t}\quad\ \text{for }\,|x|\geq r_3,
$$
then 
\begin{equation}\label{2.2b}h_{2,t}(x)\leq c_{25}e^t|x|^{-\gg_\infty-\gb_-}\quad\text{for }\,|x|\geq r_3,
\end{equation}
which implies
$$\lim_{t\to-\infty}\int_{\R^2}h_{2,t}(x) dx=0.
$$

 We claim that there exists $t_i\in\R$ such that
 \begin{equation}\label{2.1+1}
  \int_{\R^2} h_{i,t_i} (x) dx=\int_{\R^2} g_\beta (x) dx=2\pi[2(N-M)+ \beta].
 \end{equation}
From the definition of $F_i$, (\ref{2.1-1}) and the assumption on $g_\gb$,
$$\lim_{t\to+\infty}\int_{\R^2}h_{1,t}(x) dx= \infty\quad{\rm and} \ \ \lim_{t\to+\infty}\int_{\R^2}h_{2,t}(x) dx= \int_{\R^2}Wdx>\int_{\R^2} g_\beta\, dx. $$
 Since $t\mapsto \int_{\R^2}h_{i,t}(x) dx$ is continuous and increasing, it follows by the mean value theorem that there exists $t_i\in\R$ such that
$$\int_{\R^2} h_{i,t_i} (x) dx= \int_{\R^2} g_\beta(x) dx.$$
 \smallskip

\noindent {\it Step 2}. We use Lemma \ref{lm 1.2} to obtain some basic estimates on   $w_{0,i}=\Gamma\ast (g_\beta-h_{i,t_i})$, taking into account the fact that 
$\int_{\R^2} (g_\beta-h_{i,t_i})\, dx=0$ and
$$-\Delta w_{0,i} =g_\beta-h_{i,t_i}\quad {\rm in}\quad \R^2.$$
The function $g_\beta$ is smooth with compact support, the functions $h_{i,t_i}$ are locally integrable in $\BBR^2$ and satisfy
$$|h_{1,t_1}(x)|\leq c_{26}|x|^{-\gg_\infty-\gb}\quad{\rm and}\quad |h_{2,t_2}(x)|\leq c_{26}|x|^{-\gg_\infty-\gb_-}\quad\text{for $|x|$ large enough}.
$$
Since ({\ref{2La}}) holds, then, by Lemma \ref{lm 1.2}, the functions $w_{0,i}$ are uniformly bounded in $\R^2$, 
\begin{equation}\label{2Labis}
|w_{0,i}(x)|\leq c_{27}|x|^{\varrho_i}\quad\text{and }\,|\nabla w_{0,i}(x)|\leq 
c_{28}|x|^{\varrho_i-1}\quad\text{for $|x|$ large enough},
\end{equation}
where \begin{equation}\label{exp ----}{\small \varrho_1=-\frac{\gg_\infty+\gb-2}{\gg_\infty+\gb-1}\quad\text{ and } \ \varrho_2=-\frac{\gg_\infty+\gb_--2}{\gg_\infty+\gb_--1}.}\end{equation}
\smallskip

\noindent{\it Step 3}.   In order to apply the classical iterative method we have to construct suitable supersolutions and subsolutions for equation (\ref{eq 2.1}).\smallskip

\noindent {\it Construction of the supersolution.} Set $\overline v_i =(t_i)_++w_{0,i}+\norm{w_{0,i}}_{L^\infty(\R^2)}$, then
$$\BA {lll}
-\Gd\overline v_i+WF_i(\gl^\gb e^{\overline v_i})=g_\gb-WF_i(\gl^\gb e^{t_i})+WF_i(\gl^\gb e^{\overline v_i})\\
\phantom{-\Gd\overline v_i+WF_i(\gl^\gb e^{\overline v_i})}
\geq g_\gb,
\EA$$
since $F_i(\gl^\gb e^{\overline v_i})\geq F_i(\gl^\gb e^{t_i})$ as $\overline v_i\geq t_i$. Hence  $\overline v_i$ is a super solution of (\ref{eq 2.1}) for $i=1,2$.\smallskip

\noindent {\it Construction of the subsolution.} Set $\underline  v_i=-(t_i)_-+w_{0,i}-\norm{w_{0,i}}_{L^\infty(\R^2)}$, then 
$$ 
-\Gd\underline v_i+WF_i(\gl^\gb e^{\underline v_i})=g_\gb-WF_i(\gl^\gb e^{t_i})+WF_i(\gl^\gb e^{\underline v_i}) \leq g_\gb,
 $$
since $F_i(\gl^\gb e^{\underline v_i})\leq F_i(\gl^\gb e^{t_i})$ as $\underline v_i\leq t_i$. Hence  $\underline v_i$ is a subsolution of (\ref{eq 2.1}) for $i=1,2$.
As $\overline v_i>\underline  v_i $ in $\R^2$, by a standard iterating process, see \cite[Section 2.4.4]{Y1}, there exists a solution $v_i$ of
(\ref{eq 2.1})
such that
$$\underline  v_i\le v_i\le \overline v_i\quad{\rm in}\  \R^2.$$
Note that $v_i$ belongs to $C^2(\BBR^2\setminus \Gs_1)\cap C(\BBR^2)\cap L^{\infty}(\BBR^2)$. 

\smallskip
\noindent{\it Uniqueness:} Let $\tilde v_i$ be another solution of (\ref{eq 2.1}) and $w_i=\tilde v_i-v_i$, then
\begin{eqnarray*}
 \Delta (w_i^2) &=& 2w_i\Delta w_i+2|\nabla w_i|^2
 \\&\geq& 2w_i  \Delta w_i  \\
 &=& 2 w_i\left(F_i(\lambda^\beta e^{\tilde v_i})-F_i(\lambda^\beta e^{  v_i}) \right)
  \  \geq 0,
\end{eqnarray*}
hence $w_i^2$ is bounded and subharmonic in $\R^2$. Thus $w_i^2$ is a constant by Liouville's theorem, that is $\tilde v_i=v_i+C$. 
Then $F_i(\lambda^\beta e^{v_i})=F_i(\lambda^\beta e^{  v_i+C})$. Thus $C=0$ and uniqueness follows. We denote by $v_{\gb,i}$ this unique solution. \smallskip

\noindent{\it Step 3: asymptotic expansion}.  Now we shall employ Lemma \ref{lm 1.2} with $\Gf_i=WF_i(\lambda^\beta e^t)-  g_\beta$,
 where $g_\beta$ has compact support and (\ref{2La}) holds, thus 
$$\limsup_{|x|\to+\infty}\left(|\Phi_1(x)||x|^{\gamma_\infty+\beta}+|\Phi_2(x)||x|^{\gamma_\infty+\beta_-}\right)<+\infty.$$
Therefore  we have that 
 \begin{equation}\label{2.2-1}
 \limsup_{|x|\to+\infty}|\Gamma\ast \Gf_i(x)]|x|^{-\varrho_i} <+\infty.
 \end{equation}
The function  $w=v_{\gb,i}-\Gamma\ast \Gf_i$ is harmonic and bounded, hence it is constant by Louville theorem. Denote this constant by $C_{\gb,i}$, we deduce that 
 $$v_{\gb,i}=C_{\beta,i}+O(|x|^{\varrho_i})
 \
 \;{\rm as}\ \, |x|\to+\infty.$$
 The gradient estimates in (\ref{2.1 e}) are the consequences of Corollary \ref{gradecay},
which ends the proof.\hfill$\Box$

\begin{corollary}\label{flux1-2} Under the assumptions of Theorem \ref{teo 2.1} the unique solutions $v_{\gb,i}$  of problem (\ref{eq 2.1}) with $i=1,2$  respectively, satisfy the flux identity
 \begin{equation}\label{flux-1}
\int_{\BBR^2}WF_i(\gl^\gb e^{v_{\gb,i}}) dx=2\gp(2(N-M)+\gb).
 \end{equation}
\end{corollary}
{\bf Proof}. For any $R>0$, there holds
$$-\int_{|x|=R}\frac{\partial v_{\gb,i}}{\partial r}dS+\int_{B_R}WF_i(\gl e^{v_{\gb,i}}) dx=\int_{B_R} g_\gb dx.
$$
By (\ref{2.1 e}), 
$$\left|\int_{|x|=R}\frac{\partial v_{\gb,i}}{\partial r}dS\right|=O(|x|^{\rho_i})\quad\text{as }\; |x|\to+\infty,$$ 
where $\rho_i$ is defined in (\ref{exp ----}).
 The result follows from (\ref{gbeta}). \hfill$\Box$
\medskip

In the critical case $\beta=\beta^*:=2aN-2$ where $a>0$ and $0<aN\le 1$,   the problem related to (\ref{eq 1.1}) is the following
\begin{equation}\label{eq 2.1-1}
-\Delta u +  W F_2(\lambda^{\beta^*} \Lambda^{-2} e^u) =g_{\beta^*} \quad{\rm in}\;\;  \mathbb{R}^2,
\end{equation}
 where $g_{\beta^*}$ expressed by
\begin{equation}\label{gb}
 g_{\beta^*} =f_1-f_2+\beta^*\Delta \ln \lambda-2\Delta \ln \Lambda,
\end{equation}
is subject to the condition
$$\int_{\R^2} g_{\beta^*}\, dx=2\pi[2(N-M)+\beta^*],$$
and $W$ satisfies that \smallskip

\noindent $(\mathcal{W}_1)\;${\it
 The function $W$  is positive, locally H\"{o}lder continuous in $\R^2\setminus \Sigma_1$ and satisfies
$$W(x)\le c_{29}|x-p_j|^{-2n_ja}\quad{in}\quad B_{\sigma}(p_j) \quad{\rm and}\quad \limsup_{|x|\to\infty}\left( \left| |x|^{2aN}W(x)-2\right||x|\right)<+\infty, $$
\qquad\  where $c_{29}>0$, $n_ja<1$ with $j=1,\cdots,k$.}

\begin{teo}\label{teo 2.2}
Let $F_2(s)=\frac{s}{e^{\nu_1-\nu_2}+s}$, $g_{\beta^*}$ be defined in (\ref{gb}) with $\gb^*=2(aN-1)\leq  0$ and $W$ satisfies $(\mathcal{W}_1)$.
Assume furthermore that  $M<(1+a)N-1$ and set $\theta^*=\min\{3,2-\gb^*\}\geq 2$. 
Then   problem (\ref{eq 2.1-1})   has a unique bounded solution $v$ and there exists $C_*\in\R$ such that
\begin{equation}\label{2L00}\BA {lll}
\quad v_*(x)=C_*+O(|x|^{-\frac{\gth^*-2}{\gth^*-1}})\quad&\text{as }|x|\to+\infty\\[2mm]
|\nabla v(x)|= O(|x|^{-1-\frac{\gth^*-2}{\gth^*-1}})\quad&\text{as }|x|\to+\infty,
\EA\end{equation}
if $aN<1$, or
\begin{equation}\label{2L00*}\BA {lll}
\quad v_*(x)=C_*+O\left((\ln |x|)^{-4}\right)\quad&\text{as }|x|\to+\infty \\[2mm]
|\nabla v(x)|= O\left(|x|^{-1}(\ln |x|)^{-4}\right)\quad&\text{as }|x|\to+\infty,
\EA\end{equation}
if $aN=1$. 
\end{teo}
{\bf Proof}.  Notice that the assumptions $aN\leq 1$ and $M<(1+a)N-1$ imply $N-M>0$. Set
\begin{equation}\label{3.5}
\Lambda_0(x)=\frac1{1+|x|^2}\quad{\rm for\ any }  \ x\in\R^2,
\end{equation}
and for $t\in\R$, 
$$h_t(x)=  \frac{W \lambda^{\beta^*} \Lambda^{-2}  e^{ t\Lambda_0(x)}}{e^{\nu_1(x)-\nu_2(x)}+\lambda^{\beta^*} \Lambda^{-2} e^{t\Lambda_0(x)}}\quad {\rm for\ any }  \  x\in \R^2\setminus\Sigma,
 $$
with $h_t(x)=0$ for $x\in \Sigma_2$. The function $h_t(\cdot)$ is continuous  in $\R^2\setminus\Sigma_1$ and
$t\mapsto h_t(x)$ is   increasing for all $x\in\R^2\setminus\Sigma$. Direct computation implies the following properties:
 $$h_t(x)\to W\quad{\rm locally\ in}\ \R^2\setminus \Sigma_2 \quad {\rm as}\ t\to+\infty,$$
 and
 $$h_t(x)\to 0\quad{\rm locally\ in}\ \R^2\setminus \Sigma_1  \quad {\rm as}\ t\to-\infty.$$
Since $2aN\leq 2$, there holds
$$\int_{\R^2} W(x) dx=\infty.$$
Furthermore, there exist $\gt\in\BBR$ and $r_*>0$ such that for any $t\leq \gt$ and $|x|\geq r_*$, 
$$h_t(x)\leq c_{30}\frac{|x|^{-2}(\ln (|x|+1))^{-2}}{1+|x|^{2aN-2}(\ln (|x|+1))^{-2}},
$$
where $c_{30}>0$ depends on $\gt$. Since $2aN-2\leq 0$, it follows that for $|x|\geq r_*$,
\begin{equation}\label{2Lc}
h_t(x)\leq c_{31}|x|^{-2}(\ln (|x|+1))^{-2}.
\end{equation}
Hence, by the dominated convergence theorem, 
$$
\lim_{t\to-\infty}\int_{\R^2} h_t(x) dx=0.
$$
Using the fact that $t\mapsto \int_{\R^2}h_t(x) dx$ is increasing,  there exists $t_0\in\BBR$ such that
 \begin{equation}\label{2.1}
  \int_{\R^2} h_{t_0} (x) dx=2\pi[2(N-M)+\beta^*]=\int_{\R^2} g_{\beta^*} (x) dx.
 \end{equation}
 We claim that for some $c_{32}>0$,
  \begin{equation}\label{2.1x}
|w_0(x)|\leq c_{32}|x|^{-\frac{\gth^*-2}{\gth^*-1}}\quad\text{for $|x| $ large enough,}
 \end{equation}
 and if this holds true it will follow that $ \norm{w_0}_{L^\infty}<\infty$,  where   $w_0=\Gamma\ast (g_{\beta^*}-h_{t_0})$.\\
  Using (\ref{gb}), 
$$
g_{\gb^*}(x)=\frac{2}{|x|^2\Gl^2(x)}\quad\text{for $|x|\geq r_1$},
$$
  and
 $$\BA {lll}\displaystyle h_{t_0}(x)=\frac{W\gl^{\gb^*}\Gl^{-2}e^{t_0\Gl_0}}{1+\gl^{\gb^*}\Gl^{-2}e^{t_0\Gl_0}} =\frac{2|x|^{-2}e^{t_0\Gl_0}(1+O(|x|^{-1})}{\Gl^2+|x|^{\gb^*}e^{t_0\Gl_0}}\quad\text{as }\,|x|\to+\infty.
 \EA$$
 Therefore, we obtain that 
 \begin{equation}\label{2Le}\BA {lll}
 \displaystyle g_{\beta^*}-h_{t_0}=
 \frac{2e^{t_0\Gl_0}}{|x|^2}\left(\frac{\Gl^2e^{t_0\Gl_0}(1+O(|x|^{-1}))-\Gl^2-|x|^{\gb^*}e^{t_0\Gl_0}}{\Gl^2e^{t_0\Gl_0}(\Gl^2+|x|^{\gb^*}e^{t_0\Gl_0})}\right)
 \\[4mm]\phantom{g_{\beta^*}-h_{t_0}} \displaystyle
 = \frac{2e^{t_0\Gl_0}\left(\Gl^2(e^{t_0\Gl_0}-1)-e^{t_0\Gl_0}(|x|^{\gb^*}-\Gl^2O(|x|^{-1}))\right)}{|x|^2\Gl^2e^{t_0\Gl_0}(\Gl^2+|x|^{\gb^*}e^{t_0\Gl_0})}.
\EA  \end{equation}
 Since $\Gl_0(x)$ is defined by (\ref{3.5}), $e^{t_0\Gl_0}-1=O(|x|^{-2})$ at infinity. Noticing that  $\gb^*=0$ if $aN=1$, we conclude that 
 \begin{equation}\label{2Lf}|g_{\beta^*}-h_{t_0}|\leq c_{33}\max\left\{|x|^{-3}\Gl^{-2},|x|^{-2+\gb^*}\Gl^{-4}\right\}
 \leq c_{33}\left\{\BA {lll}|x|^{-\gth^*}\quad&\text{if }\;aN<1\\[2mm]
 |x|^{-2}(\ln|x|)^{-4}\quad&\text{if }\;aN=1.
\EA \right.
\end{equation}
Additionally, $\int_{\R^2} w_0 dx=0$.  Therefore,  from Lemmas \ref {lm 1.2} and \ref {lm 1.3}, we have that $w_0$ remains bounded on $\BBR^2$ and there exists $c_{34}>0$ such that  
 \begin{equation}\label{2Ld}
\big| w_0(x)\big|\leq  c_{34} (1+|x|)^{-\frac{\gth^*-2}{\gth^*-1}}\qquad\text{for all }\,x\in\BBR^2\quad \text{if $aN<1$}
 \end{equation}
  and
  \begin{equation}\label{2Ldf}
\big| w_0(x)\big|\leq  c_{34} \big(\ln(2+|x|)\big)^{-4}\qquad\text{for all }\,x\in\BBR^2 \quad \text{if $aN=1$}.
 \end{equation}
 \smallskip

\noindent{\it Existence.} We first construct a supersolution. Set
$$\overline v=(t_0)_+ +w_0+\norm{w_0}_{L^\infty(\R^2)}\quad{\rm in}\quad \R^2.$$
Since $\Lambda_0:\R^2\to (0,1]$, then $\overline v\geq t_0\Lambda_0$ in $\R^2$.
The function $t\mapsto \frac{ \Lambda^{-2} e^{ t}}{e^{\nu_1-\nu_2}+\Lambda^{-2} e^{t}}$ is increasing, therefore,
$$\frac{W \lambda^{\beta^*} \Lambda^{-2} e^{ \overline v}}{e^{\nu_1-\nu_2}+\lambda^{\beta^*} \Lambda^{-2} e^{ \overline v}} \geq \frac{W\lambda^{\beta^*}\Lambda^{-2} e^{ t_0\Lambda_0}}{e^{\nu_1-\nu_2}+\lambda^{\beta^*} \Lambda^{-2} e^{ t_0\Lambda_0}},$$
which implies,
$$
 -\Delta \overline v+\frac{W\lambda^{\beta^*} \Lambda^{-2} e^{ \overline v}}{e^{\nu_1-\nu_2}+\lambda^{\beta^*} \Lambda^{-2} e^{ \overline v}} -g_{\beta^*}
\geq  g_{\beta^*}-h_{t_0} +\frac{W\lambda^{\beta^*} \Lambda^{-2} e^{ t_0\Lambda_0}}{e^{\nu_1-\nu_2}+\lambda^{\beta^*} \Lambda^{-2} e^{ t_0\Lambda_0}}-g_{\beta^*}
=0,
$$
then $\overline v$ is a super solution of (\ref{eq 2.1}).

 Similarly we construct a subsolution by setting $\underline  v=(t_0)_-+w_0-\norm{w_0}_{L^\infty(\R^2)}$.
Using $\underline  v\le t_0\Lambda_0$ in $\R^2$ and by monotonicity, we have that 
$$\frac{W\lambda^{\beta^*} \Lambda^{-2} e^{ \underline  v}}{e^{\nu_1-\nu_2}+\lambda^{\beta^*} \Lambda^{-2} e^{ \underline  v}} \le \frac{W\lambda^{\beta^*} \Lambda^{-2} e^{ t_0\Lambda_0}}{e^{\nu_1-\nu_2}+\lambda^{\beta^*} \Lambda^{-2} e^{ t_0\Lambda_0}},$$
thus,
$$
 -\Delta \underline  v +\frac{W\lambda^{\beta^*} \Lambda^{-2} e^{ \underline  v}}{e^{\nu_1-\nu_2}+\lambda^{\beta^*} \Lambda^{-2} e^{ \underline  v}} -g_{\beta^*}\leq g_{\beta^*}-h_{t_0} +\frac{W\lambda^{\beta^*} \Lambda^{-2}  e^{ t_0\Lambda_0}}{e^{\nu_1-\nu_2}+\lambda^{\beta^*} \Lambda^{-2} e^{ t_0\Lambda_0}}-g_{\beta^*}\nonumber
=0,
$$
thus $\underline  v$ is a subsolution. 

  Since $\overline v>\underline  v $, the standard iterative process, yields the existence of a solution $v_*$ of
(\ref{eq 2.1})
such that
$$\underline  v\le v_*\le \overline v\quad{\rm in}\quad \R^2.$$
As in the proof of Theorem \ref{teo 2.1} the solutions are unique in the class of bounded solutions, a class to which $v_*$ belongs.  Put $\Gf_*=g_{\gb^*}-WF_2(\gl^*\Gl^{-2}e^{v_*})$, then $w_*=v_*-\Gg\ast \Gf_*$ is harmonic and bounded, hence it is a constant, say $C_*$. Since $ \Gf_*$ satisfies the same estimate (\ref{2Ld}), with  possibly  another constant, we deduce from Lemma \ref{lm 1.2} that for all $x\in\BBR^2$
 \begin{equation}\label{2Ld+}
 \abs {\Gg\ast \Gf_*}\leq c_{35}\left\{\BA {lll}(1+|x|)^{-\frac{\gth^*-2}{\gth^*-1}}&\quad\text{if }\; aN<1\\[2mm]
 (1+|x|)^{-2}(\ln (2+|x|))^{-4}&\quad\text{if }\; aN=1.
 \EA\right.
 \end{equation}
  This implies inequalities (\ref{2L00}) and (\ref{2L00*}) by Lemma \ref{lm 1.2} and Lemma \ref{lm 1.3} and Corollary \ref{gradecay} and Corollary \ref{gradecay*}.
 \hfill$\Box$\medskip

Similarly as Corollary \ref{flux1-2}, there holds,
\begin{corollary}\label{flux*} Under the assumptions of Theorem \ref{teo 2.2} the solutions $v_{\gb^*}$ satisfy
 \begin{equation}\label{flux-2}
\int_{\BBR^2}WF_2(\gl^{\gb^*}\Gl^{-2} e^{v_{\gb^*}}) dx=2\gp(2(N-M)+\gb^*).
 \end{equation}
\end{corollary}

From the existence and  uniqueness of solutions of (\ref{eq 2.1}) and  (\ref{eq 2.1-1}), it is easy to prove the following statements.

\begin{corollary}\label{cr 2.1} Under the assumptions of Theorem \ref{teo 2.1}, if $\overline w_i$ and $\underline w_i$ are respectively a bounded supersolution and a bounded subsolution of (\ref{eq 2.1}) such that $\underline w_i\leq \overline w_i$, then the standard iterative process will converge to the unique bounded solution $v_i$
of (\ref{eq 2.1}), and $\underline w_i\leq v_i\leq \overline w_i$. A similar result holds concerning equation (\ref{eq 2.1-1}) under the assumption of Theorem \ref{teo 2.1}.
\end{corollary}

\begin{corollary}\label{cr 2.2}
  Under the assumptions of Theorem \ref{teo 2.2}, 
the function $w_{*}:=\lambda^{\beta^*} \Lambda^{-2}+v_{*}$ where $v_{*}$ is the unique bounded solution of  (\ref{eq 2.1-1}) satisfies
  \begin{equation}\label{eq 2.1-2}
-\Delta w_{*} +  W F_2( e^{w_*}) =f_1-f_2 \quad{\rm in}\;\;  \mathbb{R}^2.
\end{equation}

\end{corollary}

\setcounter{equation}{0}
\section{Minimal solution}

In order to consider solutions $w$ of (\ref{eq 3.2}) with  asymptotic behavior $\beta\ln |x|+O(1)$, we look for $w$ under the form $w=\beta\ln |x|+v$
where  $v$ is a bounded function satisfying some related equation.
In particular, we look for non-topological solution $u_\beta$   of
problem (\ref{eq 1.1}) under the form
$$
 u_\beta =-\nu_1 +\nu_2 +\beta \ln \lambda+v_\beta\quad {or}\quad w_\beta=v_\beta+\beta\ln\lambda,
$$
where $\lambda$ is given by (\ref{lambda}) and $v_\beta$   is a bounded classical solution of
\begin{equation}\label{eq 3.1}
 \displaystyle  -\Delta v +\frac{ V \lambda^{\beta} e^{ v}}{(e^{\nu_1-\nu_2}+\lambda^{\beta} e^{v})^{1+a}} =g_\beta \quad{\rm in}\;\; \mathbb{R}^2,
\end{equation}
 with $V$ being defined in (\ref{V}) and where $g_\beta$ is defined in (\ref{2.5}). {\it Here and in what follows, we always assume that  $a=16\pi G$, $an_j< 1\ \, {\rm for}\ j=1,\cdots,k$ and $\mathcal{M}$ is the total magnetic flux  given in (\ref{tmf}).}

 We first consider the  non-topological solutions of type I for problem (\ref{eq 1.1}) in the subcritical case, which are solutions verifying
 $u(x)=\beta\ln |x|+O(1)$   as $|x|\to\infty$ with $\beta<0$.
It is equivalent to look for classical solutions of (\ref{eq 3.1})
 with $\beta<0$.

\begin{proposition}\label{pr 3.1}
Let   $N,M$ be positive integers verifying (\ref{NM1}) and $aN\le 1$, then
 for any    $\beta\in(-2(N-M), \,\gb^*)$, problem (\ref{eq 3.1}) has
a minimal bounded solution  $v_{\beta,min}$  such that
  \begin{equation}\label{3.4}
  \int_{\R^2}{V}\frac{ \lambda^{\beta} e^{ v_{\beta,min}}}{(e^{\nu_1-\nu_2}+\lambda^{\beta} e^{v_{\beta,min}})^{1+a}} dx =2\pi (2(N-M)+\beta).
  \end{equation}

\end{proposition}
{\bf Proof.}  {\it Step 1: construction of an approximating scheme}. 
We recall that
$$  {\bf P}= \frac{{V}}{e^{a(\nu_1-\nu_2)}}\quad{\rm in}\ \ \R^2\setminus\Sigma,$$
then
$$\lim_{x\to p_j} {\bf P}(x)|x-p_j|^{2n_ja}=A_0\prod_{i\neq j}^k|p_i-p_j|^{-2an_i},\quad \  \lim_{|x-q_j|\to0^+} {\bf P}(x)=0, $$
and
 \begin{equation}
 \lim_{|x|\to\infty} {\bf P}(x)|x|^{2aN}=A_0.
 \end{equation}
Since $aN\le 1$, there holds
$$\int_{\R^2}  {\bf P}(x)dx=\infty,$$
then that $ {\bf P}$ verifies the assumption $(\mathcal{W}_0)$ with $\tau_{p_j}=2n_ja<2$ and $\tau_\infty=2aN>2-2(N-M)$.
Theorem \ref{teo 2.1}-$(ii)$ implies that for any   $\beta\in(-2(N-M),\,2aN-2\})$,  the nonlinear elliptic problem
\begin{equation}\label{seq 0}
-\Delta v+ \frac{ {V}}{e^{a(\nu_1-\nu_2)}}\frac{\lambda^\beta e^{v}}{e^{\nu_1-\nu_2}+\lambda^\beta e^{v}}=g_\beta\quad{\rm in}\ \,\R^2,
\end{equation}
has a unique bounded solution $v_0$, which is continuous in $\R^2$,  smooth in $\R^2\setminus\Sigma$ and
$$\int_{\BBR^2}\left(\frac{ {V}}{e^{a(\nu_1-\nu_2)}}\frac{\lambda^\beta e^{v_0}}{e^{\nu_1-\nu_2}+\lambda^\beta e^{v_0}}-g_\beta\right) dx=0
$$
by the same argument as in Theorem \ref{teo 2.1}-(ii); then there exists a constant $C_{0,\gb}$  such that
\begin{equation}\label{seq 0bis}
\lim_{|x|\to+\infty}v_0(x)=C_{0,\gb}\quad\text{and }\,\,v_0(x)-C_{0,\gb}=O\left(|x|^{-\frac{2aN-\gb-2}{2aN-\gb-1}}\right)\quad\text{as }\; |x|\to+\infty.
\end{equation}
We set
$$ W_0={\bf P}\,\text{ and }\; W_1=\frac{ { V}}{(e^{ \nu_1-\nu_2}+\lambda^\beta e^{v_{0}})^a}
=\frac{ { e^{ a(\nu_1-\nu_2)}}}{(e^{ \nu_1-\nu_2}+\lambda^\beta e^{v_{0}})^a}W_0
\quad{\rm in}\ \ \R^2\setminus\Sigma.$$
The  function $ W_1$ is positive and H\"older continuous  in $\R^2\setminus\Gs_1$, and since
 \begin{equation}\label{LV1}
0 < W_1(x)\le W_0(x)\quad \forall x\in\R^2,
 \end{equation}
it satisfies $\CW_0$. Furthermore, as $N-M>0$, $v_0(x)\to 0$ as $|x|\to\infty$ and $\gb<0$ and therefore $W_1(x)=W_0(x)(1+o(1))$ as $|x|\to\infty$. Applying Theorem \ref{teo 2.1}-$(ii)$, with $\gg_\infty=2aN$, we see that there exists a unique bounded function $v_1$ satisfying
\begin{equation}\label{seq 1}
-\Delta v_1+W_1\frac{\lambda^\beta e^{v_1}}{e^{\nu_1-\nu_2}+\lambda^\beta e^{v_1}}=g_\beta\quad{\rm in}\ \,\R^2.
\end{equation}
Furthermore, $v_1(x)$ converges to some constant $C_{1,\gb}$ when $x\to+\infty$ and 
\begin{equation}\label{LV1bis}
v_1(x)=C_{1,\gb}+O\left(|x|^{-\frac{2aN-\gb-2}{2aN-\gb-1}}\right)\quad\text{as }\; |x|\to+\infty.
\end{equation}
Set $z=v_0-v_1$. Since the function $t\mapsto \frac{\lambda^\beta e^{t}}{e^{\nu_1-\nu_2}+\lambda^\beta e^t}$ is nondecreasing, it follows that
$$\BA {lll}\displaystyle -\Gd z^2_+=2z_+(W_1-W_0)\frac{\lambda^\beta e^{v_1}}{e^{\nu_1-\nu_2}+\lambda^\beta e^{v_1}}
-2W_0z_+\left(\frac{\lambda^\beta e^{v_0}}{e^{\nu_1-\nu_2}+\lambda^\beta e^{v_0}}-\frac{\lambda^\beta e^{v_1}}{e^{\nu_1-\nu_2}+\lambda^\beta e^{v_1}}\right)
-2|\nabla z_+|^2\\
\phantom{-\Gd z^2_+}
\leq 0.
\EA$$
Hence $z_+^2$ is subharmonic and bounded, it is therefore constant. Hence $(v_0-v_1)_+=C\geq 0$. If $C>0$ then $\sup\{v_0-v_1,0\}=C$, which implies that 
$v_0-v_1=C$. Replacing $v_0$ by $v_1+C$ we deduce from (\ref{seq 0}), (\ref{seq 1})
$$\frac{V}{(e^{ \nu_1-\nu_2}+\lambda^\beta e^{v_{1}+c})^a}\frac{\lambda^\beta e^{v_1}}{e^{\nu_1-\nu_2}+\lambda^\beta e^{v_1}}=\frac{ {V}}{e^{a(\nu_1-\nu_2)}}\frac{\lambda^\beta e^{v_1+c}}{e^{\nu_1-\nu_2}+\lambda^\beta e^{v_1+C}},
$$
which yields
$$e^C\left(e^{\nu_1-\nu_2}+\lambda^\beta e^{v_1}\right)=e^{a(\gn_1-\gn_2)}\left(e^{\nu_1-\nu_2}+\lambda^\beta e^{v_1+c}\right)^{1-a} \quad{\rm in}\ \,\R^2\setminus \Gs.
$$
Since $\gb<0$, we obtain $e^C=1$ by letting $|x|\to\infty$. Hence $C=0$ which implies $v_0\leq v_1$ in $\R^2$
and $C_{1,\beta}\geq C_{0,\beta}$.

  By induction, we suppose that for $n\geq 2$ we have constructed the sequence $\{v_{k}\}_{k<n}$ of bounded solutions to 
\begin{equation}\label{seq 0k}
-\Delta v_k+W_k\frac{\lambda^\beta e^{v_k}}{e^{\nu_1-\nu_2}+\lambda^\beta e^{v_k}}=g_\beta\quad{\rm in}\ \,\R^2,
\end{equation}
  where 
$$W_k=\frac{e^{a(\gn_1-\gn_2)}}{(e^{\gn_1-\gn_2}+\gl^\gb e^{v_{k-1}})^a}W_{0}.
$$
Then $0<W_k\leq W_{k-1}\leq ...\leq W_0$ and therefore $v_0\leq ...\leq v_{k-1}\leq   v_k$, and furthermore
\begin{equation}\label{LV1-k}
v_k(x)=C_{k,\gb}+O\left(|x|^{-\frac{2aN-\gb-2}{2aN-\gb-1}}\right)\quad\text{as }\; |x|\to+\infty.
\end{equation}
Then $v_n$ is the unique bounded solutions of 
\begin{equation}\label{seq 0n}
-\Delta v_n+W_n\frac{\lambda^\beta e^{v_n}}{e^{\nu_1-\nu_2}+\lambda^\beta e^{v_n}}=g_\beta\quad{\rm in}\ \,\R^2,
\end{equation}
where 
$$W_n=\frac{e^{a(\gn_1-\gn_2)}}{(e^{\gn_1-\gn_2}+\gl^\gb e^{v_{n-1}})^a}W_{0}\leq W_{n-1}=\frac{e^{a(\gn_1-\gn_2)}}{(e^{\gn_1-\gn_2}+\gl^\gb e^{v_{n-2}})^a}W_{0},
$$
since $v_{n-2}\leq v_{n-1}$ by induction. Furthermore, by Lemma \ref {lm 1.2} (since $\gb<0$ and $N-M>0$), 
\begin{equation}\label{LV1-n}
v_n(x)=C_{n,\gb}+O\left(|x|^{-\frac{2aN-\gb-2}{2aN-\gb-1}}\right)\quad\text{as }\; |x|\to+\infty.
\end{equation}
As above the function $(v_{n-1}-v_n)_+^2$ is subharmonic and bounded, hence it is constant, which implies $v_{n-1}=v_n+C$, $C\geq 0$. Then from the equations satisfied by $v_n$ and $v_{n-1}$, 
$$1\geq 
\frac{W_n}{W_{n-1}}=e^C\frac{e^{\gn_1-\gn_2}+\gl^\gb e^{v_{n}}}{e^{\gn_1-\gn_2}+\gl^\gb e^Ce^{v_{n}}}\geq 1.$$
Hence $e^C=1$, then $C=0$ and $v_{n-1}\leq v_n$. Consequently $n\mapsto C_{n,\beta}$ is increasing. \smallskip

  Let $R>1$ be such that supp $\!\!(g_\gb)\subset B_R$ and $\Gth(x):=\Gg\ast |g_\beta|(x)=\frac{1}{2\gp}\int_{B_R}|g_\gb(x)|\ln|x-y|dy
$. For $|x|\geq R+1$, one has $1\leq |x-y|\leq |x|+R$, hence $$0\leq \ln|x-y|\leq\ln (|x|+R)\leq \ln |x|+\frac{R}{|x|}\leq \ln |x|+1,$$ therefore
$$0\leq \Gth(x)\leq \frac{R^2\norm{g_\gb}_{L^\infty}}{2}\left(\ln|x|+1\right)\quad\text{if }\;|x|\geq R+1.
$$ 
Since $|\Gth|$ is bounded from above on $B_{R+1}$ by some constant $c_{36}$, we deduce
\begin{equation}\label{LV3}
|\Gth(x)|\leq c_{36}+\frac{R^2\norm{g_\gb}_{L^\infty}}{2}\left(\ln_+ |x|+1\right)\quad\text{for all }x\in\BBR^2.
\end{equation}
Set $z=v_n-\Gth$, then 
$$-\Gd z^2_+\leq -2z_+\Gd z=-2z_+W_n\frac{\lambda^\beta e^{v_n}}{e^{\nu_1-\nu_2}+\lambda^\beta e^{v_n}}\leq 0.
$$
The function $z_+$ has compact support because of (\ref{LV3}). It is subharmonic, nonnegative and bounded, hence it is constant with zero value necessarily, hence, for any $n\in\BBN$, 
\begin{equation}\label{LV4}
v_0(x)\leq v_n(x)\leq c_{36}+\frac{R^2\norm{g_\gb}_{L^\infty}}{2}\left(\ln_+ |x|+1\right)\quad\text{for all }x\in\BBR^2.
\end{equation}
For $\ge>0$ set 
$$w_\ge(x)=\ge\ln |x|+c_{36}+\frac{R^2\norm{g_\gb}_{L^\infty}}{2}\left(\ln R+1\right).
$$
Then $w_\ge$ is harmonic in $\overline {B_R}^c$. It is larger than $v_n$ for $|x|=R$ and also at infinity, since $v_n$ is bounded. If we set $Z=v_n-w_\ge$, then as above the function 
$Z_+^2$ is subharmonic, nonnegative and bounded in $B_R^c$. Since it vanishes for $|x|=R$, its extension $\gz$ by $0$ in $\overline {B_R}$ is still subharmonic nonnegative and bounded. It is therefore constant. Since it vanishes at infinity, it is identically $0$. Hence $v_n-w_\ge\leq 0$. Letting $\ge\to 0$ 
we obtain
\begin{equation}\label{LV4x}
v_0(x)\leq v_n(x)\leq c_{36}+\frac{R^2\norm{g_\gb}_{L^\infty}}{2}\left(\ln R+1\right)\quad\text{for all }x\in\BBR^2.
\end{equation}
Combining Lemma \ref{lm 1.2} with (\ref{LV4x}) for $n=1$ we infer
\begin{equation}\label{LV4xx}
 \abs{v_n(x)-C_{n,\beta}}\leq c_{37}(1+|x|)^{-\frac{2aN-\gb-2}{2aN-\gb-1}}\quad\text{for all }\; x\in\BBR^2,
\end{equation}
where $c_{37}>0$ is independent of $n$. By Lemma \ref{lm 1.2}
\begin{equation}\label{LV6-1}\abs{\nabla v_n(x)}\leq c_{38}|x|^{-1-\frac{2aN-\gb-2}{2aN-\gb-1}}\quad\text{for $|x|$ large enough}, 
\end{equation}
then, 
$$-\int_{|x|=R}\frac{\partial v_n}{\partial r}dS+\int_{B_R}W_n\frac{\lambda^\beta e^{v_n}}{e^{\nu_1-\nu_2}+\lambda^\beta e^{v_n}}
dx=\int_{B_R}g_\gb dx.
$$
By (\ref{LV6-1}), the first integral tends to $0$ when $R\to+\infty$, therefore
\begin{equation}\label{LV7}
\int_{\BBR^2}W_n\frac{\lambda^\beta e^{v_n}}{e^{\nu_1-\nu_2}+\lambda^\beta e^{v_n}}
dx=\int_{\BBR^2}g_\gb dx=2\gp(2(N-M)+\gb).
\end{equation}
Set $\displaystyle v_{\beta,min}=\lim_{n\to\infty}v_n$ and $\displaystyle  C_\beta=\lim_{n\to\infty}C_{n,\beta}$, then
$$W_n\to W_\infty:=\frac{e^{a(\gn_1-\gn_2)}}{(e^{\gn_1-\gn_2}+\gl^\gb e^{v_{\beta,min}})^a}W_{0}
$$
and
\begin{equation}\label{LV5}
 \abs{v_{\beta,min}(x)-C_\beta}\leq c_{39}(1+|x|)^{-\frac{2aN+2(N-M)-\gb-2}{2aN+2(N-M)-\gb-1}}\quad\text{for all }\; x\in\BBR^2.
\end{equation}
Furthermore
$$0\leq W_n\frac{\lambda^\beta e^{v_n}}{e^{\nu_1-\nu_2}+\lambda^\beta e^{v_n}}
\leq W_0\frac{\lambda^\beta e^{v_{\beta,min}}}{e^{\nu_1-\nu_2}+\lambda^\beta e^{v_{\beta,min}}}.
$$
The right-hand side of the above inequality is an integrable function, therefore 
$$W_n\frac{\lambda^\beta e^{v_n}}{e^{\nu_1-\nu_2}+\lambda^\beta e^{v_n}}\to W_0\frac{ e^{a(\gn_1-\gn_2)} \lambda^{\beta} e^{ v_{\beta,min}}}{(e^{\nu_1-\nu_2}+\lambda^{\beta} e^{v_{\beta,min}})^{1+a}}\quad\text{in }\,L^1(\BBR^2)\quad\text{as }\;n\to+\infty.$$
This implies that $v_{\beta,min}$ is a weak solution of (\ref{eq 1.1}) and relation (\ref{3.4}) holds.\smallskip

\noindent{\it Step 2: $v_{\beta,min}$ is minimal among the bounded solutions}. Let $\tilde v$ be any bounded solution. Then 
$\frac{V}{(e^{\gn_1-\gn_2}+\gl^\gb e^{\tilde v})^a}\leq \frac{V}{e^{a(\gn_1-\gn_2)}}$, and by uniqueness, it implies $v_0\leq\tilde v$. Hence 
$\frac{V}{(e^{\gn_1-\gn_2}+\gl^\gb e^{\tilde v})^a}\leq \frac{V}{(e^{\gn_1-\gn_2}+\gl^\gb e^{v_0})^a}$ and therefore $v_1\leq \tilde v$. By induction we obtain  $v_n\leq \tilde v$ and finally $v_{\beta,min}\leq \tilde v$. \smallskip

\noindent{\it Step 3: asymptotic behaviour}. Put 
$$F=g_\gb-\frac{V\gl^\gb e^{v_{\beta,min}}}{(e^{\gn_1-\gn_2}+\gl^\gb e^{v_{\beta,min}})^{1+a}}.
$$
Then 
$\int_{\BBR^2}F\,dx=0$ and $|F(x)|\leq c_{40}|x|^{-(2aN-\gb)}\quad\text{for }\; |x|\geq r_0.$ So we have (\ref{3.4}), and applying   Lemma \ref{lm 1.2} yields the estimate  
$$|\Gg\ast F|\leq c_{41}(1+|x|)^{-\frac{2aN-\gb-2}{2aN-\gb-1}}\quad\text{for all }\; x\in\BBR^2.
$$
Therefore $w=v_{\beta,min}-\Gg\ast F$ is harmonic and bounded in $\BBR^2$. It is therefore constant. This implies
\begin{equation}\label{LV6}
v_{\beta,min}=C_\gb+O(|x|^{-\frac{2aN-\gb-2}{2aN-\gb-1}})\quad\text{as }\,|x|\to+\infty.
\end{equation}
This ends the proof.\hfill$\Box$\medskip

\noindent{\bf  Proof of Theorem \ref{teo 0} part $(i)$.}   Let
$$
 u_{\beta,min} =-\nu_1 +\nu_2 +\beta \ln \lambda+v_{\beta,min},
$$
where $v_{\beta,min}$ is the minimal bounded solution of (\ref{eq 3.1}) obtained in Proposition \ref{pr 3.1}. Then $ u_{\beta,min}$ is the minimal   non-topological solution of type I of (\ref{eq 1.1})
in the sense that
$$u_{\beta,min}(x)-\beta\ln|x|=O(1)\quad{\rm as}\quad |x|\to+\infty.$$    Moreover, $ u_{\beta,min}$  verifies   (\ref{1.3}) and its total magnetic flux is $2\pi[2(N-M)+\beta]$
by (\ref{3.4}).
  \hfill$\Box$


\setcounter{equation}{0}
\section{Critical-minimal solutions}

\subsection{Non-topological solutions }

If $N>M$ and $aN<1$, we recall that by Theorem \ref{teo 2.1}, for any $\gb\in (-2(N-M),\gb^*)$, with $\gb^*=2(aN-1)<0$, there exists a unique bounded solution $v_\gb$ to equation
\begin{equation}\label{3.2-0}
- \Delta v +\frac{V}{e^{a(\nu_1-\nu_2)}}\frac{\lambda^\beta e^v}{e^{ \nu_1-\nu_2}+  \lambda^\beta e^{v} } =g_\beta\quad{\rm in}\;\; \mathbb{R}^2,
\end{equation}
and by Theorem \ref{teo 2.2}, there exists a unique bounded solution $v_{\gb^*}$ to 
\begin{equation}\label{3.2-0^*}
- \Delta v +\frac{V}{e^{a(\nu_1-\nu_2)}}\frac{\lambda^{\beta^*} \Lambda^{-2} e^v}{e^{ \nu_1-\nu_2}+  \lambda^{\beta^*} \Lambda^{-2} e^{v} } =g_{\gb^*}\quad{\rm in}\;\; \mathbb{R}^2.
\end{equation}
For $\beta\in (-2(N-M),\ \beta^*)$, we first set
 \begin{equation}\label{w b-}
 w_{\beta}=v_{\beta}+\beta \ln \lambda.
\end{equation}
Then $w_{\beta}$ is a the unique solution of
\begin{equation}\label{3.2-}
- \Delta w +\frac{V}{e^{a(\nu_1-\nu_2)}}\frac{e^w}{e^{ \nu_1-\nu_2}+  e^{w} } =f_1-f_2\quad{\rm in}\;\; \mathbb{R}^2
\end{equation}
such that $w-\gb\ln\gl$ is bounded in $\mathbb{R}^2$. \smallskip

\noindent  When $\gb={\gb^*}$, we set 
\begin{equation}\label{w b--}
 w_{\beta^*}=v_{\beta^*}+\beta^* \ln \lambda-2\ln\Lambda,
\end{equation}
thus $w_{\beta^*}$ is the unique solution of (\ref{3.2-}) 
such that $w-\gb^*\ln\gl-2\ln\Gl$ remains bounded in $\mathbb{R}^2$.

\begin{proposition}\label{pr 2.2}
Under the assumptions of Theorem \ref{teo 2.2}, the mapping $\beta\mapsto w_\beta$
is increasing for $\gb\in (-2(N-M),\ \beta^*)$ and
$$w_{\beta^*}=\sup\left\{ w_\beta\ {\rm in}\ \R^2:\ \beta\in (-2(N-M),\ \beta^*)\right\}.$$

\end{proposition}
{\bf Proof.} If $\beta^*>\gb>\gb'>-2(N-M)$, the function $z=w_{\gb'}-w_{\gb}$ is negative in $B_R^c$ for some $R>0$. Hence
$$-\Gd z^2_+=-2z_+\Gd z-2|\nabla z_+|^2\leq -\frac{2V}{e^{a(\nu_1-\nu_2)}}\left(\frac{e^w_{\gb'}}{e^{ \nu_1-\nu_2}+  e^{w_{\gb'}} }-
\frac{e^w_{\gb}}{e^{ \nu_1-\nu_2}+  e^{w_{\gb}} }\right)(w_{\gb'}-w_{\gb})\leq 0.
$$
Hence $z^2_+$ is a nonnegative and bounded subharmonic function in $\mathbb{R}^2$, it is therefore constant. Since it vanishes in 
$B_R^c$, it is identically $0$, which yields $w_{\gb'}\leq w_{\gb}$. Actually the inequality is strict since it is the case at infinity and there cannot exist 
$x_0\in \BBR^2$ such that $w_{\gb'}(x_0)= w_{\gb}(x_0)$, because of the strong maximum principle. Similarly, if $\gb<\gb^*$, there holds by (\ref{w b-}) and 
(\ref{w b--}),
$$(w_{\gb}-w_{\gb^*})(x)=(\gb-\gb^*)\ln |x|+2\ln(\ln |x|)+O(1)\quad \text{as }\,|x|\to+\infty.$$ 
Hence $z^2_+=(w_{\gb}-w_{\gb^*})^2_+$ is subharmonic nonnegative and bounded, hence it is constant and necessarily with value zero. Therefore 
$w_{\gb}\leq w_{\gb^*}$, and 
actually $w_{\gb}<w_{\gb^*}$ by the strong maximum principle. We set
$$\tilde w_{\beta^*}:=\sup\{w_\beta: \beta\in (-2(N-M),\beta^*) \}=\lim_{\gb\uparrow\gb^*}w_\beta.
$$
Then $\tilde w_{\beta^*}\leq w_{\beta^*}$ and $w_{\beta^*}$ is a solution of (\ref{3.2-}). By the strong maximum principle, either 
$\tilde w_{\beta^*}< w_{\beta^*}$ or $\tilde w_{\beta^*}= w_{\beta^*}$. In order to identify $w_{\beta^*}$, we use the flux identities obtained 
in Corollaries \ref{flux1-2} and \ref{flux*}, replacing $v_{\beta, 2}$ and $v_{\beta^*}$ by their respective expressions from (\ref{w b-}) and (\ref{w b--}):
$$\CM(w_\gb)=\int_{\BBR^2}\frac{V}{e^{\gn_1-\gn_2}}\frac{ e^{w_\gb}}{e^{\gn_1-\gn_2}+ e^{w_\gb}}=2\gp(2(N-M)+\gb),
$$
and 
$$\CM(w_{\gb^*})=\int_{\BBR^2}\frac{V}{e^{\gn_1-\gn_2}}\frac{ e^{w_{\gb^*}}}{e^{\gn_1-\gn_2}+e^{w_{\gb^*}} }dx=2\gp(2(N-M)+\gb^*).$$
Since the mapping $\beta\mapsto \frac{e^{w_\beta}}{1+e^{w_\beta}}$ is   increasing, there holds by the monotone convergence theorem,
\begin{equation}\label{bof}\BA {lll}\displaystyle
\mathcal{M}(\tilde w_{\beta^*})=\int_{\R^2}\frac{V}{e^{a(\nu_1-\nu_2)}}\frac{e^{\tilde w_{\beta^*}}}{e^{\nu_1-\nu_2}+e^{\tilde w_{\beta^*}}}dx
\\[4mm]\displaystyle
\phantom{\mathcal{M}(\tilde w_{\beta^*})}
=\lim_{\beta\uparrow \beta^*} \int_{\R^2}\frac{V}{e^{a(\nu_1-\nu_2)}}\frac{e^{ w_{\beta}}}{e^{\nu_1-\nu_2}+e^{ w_{\beta}}}dx
\\[4mm]\phantom{\mathcal{M}(\tilde w_{\beta^*})}
=2\pi(2(N-M)+\beta^*)\\[2mm]\phantom{\mathcal{M}(\tilde w_{\beta^*})}\displaystyle
=\int_{\R^2}\frac{V}{e^{a(\nu_1-\nu_2)}}\frac{e^{w_{\beta^*}}}{e^{\nu_1-\nu_2}+e^{\tilde w_{\beta^*}}}dx,
\EA\end{equation}
where $\frac{e^{\tilde w_{\beta^*}}}{e^{\nu_1-\nu_2}+e^{\tilde w_{\beta^*}}}\leq \frac{e^{w_{\beta^*}}}{e^{\nu_1-\nu_2}+e^{w_{\beta^*}}}$.  This implies that
$\tilde w_{\beta^*}= w_{\beta^*}$ almost everywhere and actually everywhere by continuity.\qeda

\subsection{Proof of Theorem \ref{teo 0-c}}

If $v_{\beta,min}$ is the minimal bounded solution of (\ref{eq 3.1}) obtained in Proposition \ref{pr 3.1}, we set
 \begin{equation}\label{w b}
 w_{\beta,min}=v_{\beta,min}+\beta \ln \lambda\quad{\rm  in}\quad \R^2.
\end{equation}
Then $w_{\beta,min}$ is a solution of
\begin{equation}\label{3.2+}
\left\{\BA{lll}
\displaystyle - \Delta w +\frac{ V e^{w  }}{(e^{\nu_1-\nu_2}+  e^{w})^{1+a}} =f_1-f_2\quad{\rm in}\;\; \mathbb{R}^2,
\\[3.5mm]\phantom{ -- }
 w =\beta\ln\lambda+O(1)\qquad{\rm as}\  |x|\to+\infty.
\EA\right.
\end{equation}
Since $v_{\beta,min}$ is the minimal bounded solution of (\ref{eq 3.1}),  $w_{\beta,min}$ is the minimal solution of (\ref{3.2+}). 

Furthermore, $v_{\beta,min}$ is the limit 
of the increasing sequence  of  the bounded solutions $\{v_n\}$ of (\ref{seq 0n}), therefore  $w_{\beta,min}$  is the limit 
of the increasing sequence $\{w_{\gb,n}\}:=\{v_{n}+\beta \ln \lambda\}$ of the solutions   of
\begin{equation} \label{3.1-n}
\left\{\begin{array}{lll}
\displaystyle -\Delta w_{\gb,n} + \frac{V}{(e^{ \nu_1-\nu_2}+  e^{w_{\gb,n-1}})^a}\frac{  e^{w_{\gb,n}}}{e^{ \nu_1-\nu_2}+  e^{w_{\gb,n}} }=f_1-f_2\quad{\rm in}\;\; \mathbb{R}^2,
\\[3.5mm]\phantom{ --}
 w_{\gb,n} =\beta\ln\lambda+O(1)\qquad {\rm as}\  |x|\to\infty.
 \end{array}\right.
\end{equation}
By the comparison principle, the mapping $\beta\in(-2(N-M),\, \beta^*)\mapsto w_{\beta,n}$ is increasing for any $n$, 
and this is also true for $\beta\mapsto w_{\beta,min}$. By (\ref{LV4x}) there holds for any $n\in\BBN$, 
\begin{equation} \label{3.1-x}
w_0(x)\leq w_{\beta,n}(x)\leq w_{\beta,min}(x)\leq c_{42}+\frac{R^2\norm{g_\gb}_{L^\infty}}{2}(\ln R+1)+\gb\ln\gl(x)\qquad{\rm in}\;\; \mathbb{R}^2.
\end{equation}
 \smallskip

\noindent {\it Uniformly  upper bound for $\{w_{\beta,min}\}_\beta$.}  Let $\overline v_2=\Gamma\ast  g_{\beta^*} $,
 then  $\overline v_2 =\Gamma\ast(f_1-f_2)+\beta^*\ln \lambda-2\ln\Lambda$ and
 $$\lim_{|x|\to\infty}\frac{\overline v_2(x)}{  \ln|x|}=2(N-M)+\beta^*>0.$$
Since $\overline v_2$ is a super solution of (\ref{3.1-n}), we have by comparison
$$w_{\beta,n}\le \overline v_2\quad{\rm in}\quad \R^2,$$
which implies that for any $\beta\in (-2(N-M),\, \beta^*)$
$$w_{\beta,min}\le \overline v_2\quad{\rm in}\quad \R^2.$$
  Hence there exists  $\displaystyle w_{\beta^*,min}=\lim_{\gb\uparrow\gb^*}w_{\beta,min}$ and 
$$w_{\beta^*,min}\le \overline v_2\quad{\rm in}\quad \R^2,$$
and therefore
\begin{equation}\label{3.9}
w_{\beta^*,min}(x)\le \beta^*\ln |x|-2\ln\ln |x|+C
\end{equation}
for some $C\in\R$.\smallskip

\noindent {\it Lower bound for $w_{\beta^*,min}$.} From Proposition \ref{pr 2.2}, the equation
\begin{equation}
-\Delta w + \frac{V}{e^{ a(\nu_1-\nu_2)} }\frac{  e^{w}}{e^{ \nu_1-\nu_2}+  e^{w} }=f_1-f_2\quad{\rm in}\ \,\R^2
\end{equation}
 has a unique solution $w_{\beta^*}$, with the following asymptotic behavior
 $$w_{\beta^*}(x)=\beta^*\ln |x|-2\ln\ln|x|+O(1)\quad{\rm as}\quad |x|\to+\infty,$$
and $w_{\beta^*}$  is the limit of the solutions $w_{\beta}$ of (\ref{3.2-}) for $\beta\in(-2(N-M), \beta^*)$ satisfying
 $$w_{\beta}(x)=\beta\ln |x| +O(1)\quad{\rm as}\ \, |x|\to+\infty.$$
Since $w_{\beta}$ is a subsolution for (\ref{3.2+}) it is bounded from above by $w_{\beta,min}$ by the same comparison method as the ones used previously. Therefore $w_{\beta^*}\le w_{\beta^*,min}.$ Combining  (\ref{3.9}) with the expression of $w_{\beta^*}$ given in (\ref{3.9}), we infer that 
$$w_{\beta^*,m}(x)= \beta^*\ln |x|-2\ln\ln |x|+O(1)\quad {\rm as}\ |x|\to+\infty. $$
Clearly the flux identity holds as in the previous theorem, which ends the proof.\hfill$\Box$

\setcounter{equation}{0}
\section{   Multiple  solutions  }

\subsection{Non-topological solutions}

Let $\beta\neq 0$ and $u_\beta$ be a   solution   of
problem (\ref{eq 1.1}) with the asymptotic behavior
$$u_\beta(x)=\beta \ln |x|+O(1)\quad\text{as }\;|x|\to+\infty.$$
Then $u_\beta$ can written under the form
$$
 u_\beta =-\nu_1 +\nu_2 +\beta \ln \lambda+v_\beta,
$$
where  $v_\beta$ is a bounded solution of the following equation equivalent to (\ref{eq 3.1})
\begin{equation}\label{eq 4.1}
 \displaystyle  -\Delta v +W_\beta\frac{   e^{ v}}{(e^{\nu_1-\nu_2}\lambda^{-\beta}+ e^{v})^{1+a}} =g_\beta \quad{\rm in}\;\; \mathbb{R}^2
\end{equation}
with $W_\beta=V\lambda^{-a\beta}$, and where $g_\beta$ is expressed by
$$
 g_\beta =f_1-f_2+\beta\Delta \ln \lambda.
$$
Note that it is a smooth function with compact support in $B_{r_0}(0)$ and it verifies
$$\int_{\R^2} g_\beta\, dx=2\pi[2(N-M)+ \beta].$$
As for $W_\beta$ it satisfies
$$\lim_{x\to p_j}W_\beta(x)= A_0(\prod_{i\neq j}|p_j-p_i|^{2n_i} )^{-a}\,,\quad \lim_{x\to q_j}W_\beta(x)= 0\;\text { and }\;\lim_{|x|\to\infty} W_\beta(x) |x|^{2aN+a\beta }= A_0.$$

The existence of multiple solutions  states as follows:

\begin{proposition}\label{pr 4.1}
Let   $N,\,M$ be positive integers and $ \gb^{\sharp}$ be given in (\ref{1.1+q}). 
Then for any $\gb>\gb^{\sharp}$
problem (\ref{eq 4.1}) possesses a sequence of solutions  $v_{\beta,i}$ such that
  \begin{equation}\label{4.1}
  \int_{\R^2}W_\beta\frac{   e^{ v_{\beta,i}}}{(e^{\nu_1-\nu_2}\lambda^{-\beta}+ e^{v_{\beta,i}})^{1+a}} dx =2\pi[2(N-M)+\beta],
  \end{equation}
  and
 $$v_{\beta,i}(x) = C_{\gb,i}+O(|x|^{-\frac{a\gb+2aN-2}{a\gb+2aN-1}})\quad\text{as }\;|x|\to+\infty$$
 with
 $$C_{\gb,i}\to\infty\ \ {\rm as}\quad i\to+\infty.$$

\end{proposition}
{\bf Proof.} By Theorem \ref{teo 2.1}, for any $A>0$ the equation
\begin{equation}\label{eq 4.2}
-\Delta w+   e^{-A(1+a)} W_\beta    e^{ w}  =g_\beta\quad{\rm in}\ \,\R^2
\end{equation}
has a unique bounded solution $w_A$.
We note that
$$w_A=w_0+A(1+a),$$
where $w_0$ is the bounded solution of (\ref{eq 4.2}) with $A=0$.
Note that   for any $A\geq A^*=a^{-1}\norm{w_0}_{L^\infty(\R^2)}$,
$$w_A\geq A\quad{\rm in}\quad \R^2.$$
{\it Step 1: construction of an approximating sequence}. We set $v_0:=w_A$ and  define $H_t(.)$ by
$$H_0(t,\cdot)=\arraycolsep=1pt\left\{
\begin{array}{lll}
 \displaystyle  A_0\qquad \quad &\text{in }\, \Sigma_1,
\\[1.5mm]\phantom{ }
 \displaystyle 0\qquad  \quad   &\text{in }\,\Sigma_2,
\\[2mm]\phantom{ }
 \displaystyle W_\beta\frac{   e^{ t}}{(e^{\nu_1-\nu_2}\lambda^{-\beta}+ e^{v_{0}})^{1+a}} \quad  &\text{in }\, \R^2\setminus\Sigma.
\end{array}
\right.$$
Under the assumptions, $H_0(t,\cdot)\in L^\gd(\BBR^2)$ for some $\gd>1$ and there exists a unique (and explicit) real number $t_1$ such that 
$$\int_{\BBR^2}H_0(t_1,x)dx=2\gp(2(N-M)+\gb).
$$
We construct first a bounded solution $v_1$ of 
\begin{equation}\label{eq 4.3n=1}
-\Delta v+ W_\beta\frac{   e^{ v}}{(e^{\nu_1-\nu_2}\lambda^{-\beta}+ e^{v_{0}})^{1+a}} =g_\beta\quad{\rm in}\ \,\R^2.
\end{equation}
 We set 
$$w_1=\Gf\ast(g_\gb-H_0(t_1,\cdot)).
$$
By Lemma \ref{lm 1.2}, $w_1$ is bounded. Put $\overline v=\norm{v_0}_{L^\infty}+w_1+\norm{w_1}_{L^\infty}+|t_1|$. Then
$$\frac{e^{\overline v}}{(e^{\nu_1-\nu_2}\lambda^{-\beta}+ e^{v_{0}})^{1+a}}\geq \frac{e^{t_1}}{(e^{\nu_1-\nu_2}\lambda^{-\beta}+ e^{v_{0}})^{1+a}},
$$
therefore
$$\BA {lll}\displaystyle 
-\Gd \overline v+W_\gb\frac{e^{\overline v}}{(e^{\nu_1-\nu_2}\lambda^{-\beta}+ e^{v_{0}})^{1+a}}-g_\gb\geq
g_\gb-H(t_1,\cdot)+W_\gb\frac{e^{t_1}}{(e^{\nu_1-\nu_2}\lambda^{-\beta}+ e^{v_{0}})^{1+a}}-g_\gb\geq 0.
\EA$$
Hence $\overline v$ is a supersolution of (\ref{eq 4.3n=1}). Since 
$$-\Gd v_0+W_\gb\frac{e^{v_0}}{(e^{\nu_1-\nu_2}\lambda^{-\beta}+ e^{v_{0}})^{1+a}}-g_\gb\leq 
W_\gb e^{v_0}\left(\frac{1}{(e^{\nu_1-\nu_2}\lambda^{-\beta}+ e^{v_{0}})^{1+a}}-\frac{1}{e^{(1+a)A}}\right)
\leq 0,$$
$v_0$ is a subsolution of (\ref{eq 4.3n=1}) dominated by $\overline v$. Hence there exists a solution $v=v_1$ of (\ref{eq 4.3n=1}) satisfying 
$$v_0\leq v_1\leq \overline v.
$$
Since $a\gb+2aN>2$, we have from Lemma \ref{lm 1.2}
$$v_1(x)=C_{1,\gb}+O(|x|^{-\frac{a\gb+2aN-2}{a\gb+2aN-1}})\quad\text{as }\;|x|\to\infty.
$$ 
We define a sequence $\{v_n\}_{n\in\BBN}$ with $v_0=w_A$  and $v=v_n$ is the bounded solution of
\begin{equation}\label{eq 4.3}
-\Delta v+ W_\beta\frac{   e^{ v}}{(e^{\nu_1-\nu_2}\lambda^{-\beta}+ e^{v_{n-1}})^{1+a}} =g_\beta\quad{\rm in}\ \,\R^2.
\end{equation}
Assume that we have proved the existence and boundedness of the functions $v_k$ for $k<n$ and that  there holds $v_0\leq v_1\leq ...\leq v_{n-1}$. We define 
$H_{n-1}(t,.)$ by 
$$H_{n-1}(t,.)=\arraycolsep=1pt\left\{
\begin{array}{lll}
 \displaystyle  A_0\qquad \quad &\text{in }\, \Sigma_1
\\[2mm]\phantom{ }
 \displaystyle 0\qquad  \quad   &\text{in }\,\Sigma_2
\\[2.5mm]\phantom{ }
 \displaystyle W_\beta\frac{   e^{ t}}{(e^{\nu_1-\nu_2}\lambda^{-\beta}+ e^{v_{n-1}})^{1+a}} \quad  &\text{in }\, \R^2\setminus\Sigma,
\end{array}
\right.$$
and denote by $t_n$ the unique real number such that
$$\int_{\BBR^2}H_{n-1}(t,x)dx=2\gp(2(N-M)+\gb).
$$
Since $v_0\leq v_1\leq ...\leq v_{n-1}$, there holds $t_0<t_1<...<t_n$. If we set $w_n=\Gg\ast(g_\gb-H_{n-1}(t_n,.))$, clearly 
 $\overline v_n:=\norm{v_{n-1}}_{L^\infty}+w_n+\norm{w_n}_{L^\infty}+|t_n|$ is a supersolution. Furthermore
 $$\BA {lll}\displaystyle-\Gd v_{n-1}+W_\gb\frac{e^{v_{n-1}}}{(e^{\nu_1-\nu_2}\lambda^{-\beta}+ e^{v_{{n-1}}})^{1+a}}-g_\gb\\[4mm]
 \phantom{---------}\leq \displaystyle
W_\gb e^{v_{n-1}}\left(\frac{1}{(e^{\nu_1-\nu_2}\lambda^{-\beta}+ e^{v_{n-1}})^{1+a}}-\frac{1}{(e^{\nu_1-\nu_2}\lambda^{-\beta}+ e^{v_{n-2}})^{1+a}}\right)
\leq 0.
\EA $$
Hence $v_{n-1}$ is a subsolution. A solution $v=v_n$ of (\ref{eq 4.3}) satisfying $v_{n-1}\leq v_n\leq \overline{v}_n$ exists. It  is bounded and satisfies 
\begin{equation}\label{eq 4.4}
v_n(x)=C_{n,\gb}+O(|x|^{-\frac{a\gb+2aN-2}{a\gb+2aN-1}})\quad\text{as }\;|x|\to\infty
\end{equation}
for some $C_{n,\gb}$. Furthermore the sequence $\{C_{n,\gb}\}$ is nondecreasing, and by Corollary \ref{gradecay}
\begin{equation}\label{eq 4.4+}
|\nabla v_n(x)|=O(|x|^{-1-\frac{a\gb+2aN-2}{a\gb+2aN-1}})\quad\text{as }\;|x|\to\infty.
\end{equation}

\noindent{\it Uniformly  upper bound for $\{v_n\}_n$.} Let $\overline v_\gb=\Gamma\ast g_\beta$, then it is a supersolution of (\ref{eq 4.3}) for any $n\in\BBN$ and satisfies
$$\lim_{|x|\to\infty}\frac{\overline v_\gb(x)}{ \ln |x|}=2(N-M)+\beta.$$
This implies that for any for any $\ge>0$, there exists $C_\ge>0$ such that 
\begin{equation}\label{eq 4.5}\overline v_\gb(x)\leq (2(N-M)+\beta+\ge)\ln( |x|+1)+C_\ge\quad\text{in }\;\BBR^2.\end{equation}
Note that $\overline v_\gb$ is a super solution of (\ref{eq 4.3}) and by the comparison principle
$$v_n\le \overline v_\gb\quad{\rm in}\quad \R^2.$$
Therefore the limit of the sequence $\{v_n\}$ as $n\to\infty$ exists. As it depends also on $A$, we denote it by $v_{\beta,A}$ and there holds
$$v_{\beta,A}\leq \overline v_\gb\quad{\rm in}\quad \R^2.$$
Furthermore $v_{\beta,A}$ is a locally bounded solution of (\ref{eq 4.1}) which satisfies
\begin{equation}\label{eq 4.6}
A\leq v_{\beta,A}(x)\leq (2(N-M)+\beta+\ge)\ln( |x|+1)+C_\ge\quad\text{in }\;\BBR^2.
\end{equation}
Because of the above lower estimate, the functions $x\mapsto \frac{   e^{ v_{n}(x)}}{(e^{\nu_1-\nu_2}\lambda^{-\beta}+ e^{v_{n}(x)})^{1+a}}$ are upper bounded on $\BBR^2$ by some constant depending on $A$ and $\gb$ but independent of $n$, and this estimate holds true if $v_n$ is replaced by $v_{\beta,A}$. Hence
for any $R>0$,
$$
-\int_{|x|=R}\frac{\partial v_n}{\partial r}dS+\int_{B_R}W_\gb\frac{ e^{v_n}}{e^{\nu_1-\nu_2}\gl^{-\gb}+ e^{v_n}}
dx=\int_{B_R}g_\gb dx.
$$
By (\ref{eq 4.4+}) the integral term on $|x|=R$ tends to $0$ when $R\to\infty$, therefore 
\begin{equation}\label{eq 4.7}
\int_{\BBR^2}W_\gb\frac{ e^{v_n}}{e^{\nu_1-\nu_2}\gl^{-\gb}+ e^{v_n}}
dx=\int_{\BBR^2}g_\gb dx=2\gp(2(N-M)+\gb).
\end{equation}
Since $\frac{ e^{v_n}}{e^{\nu_1-\nu_2}\gl^{-\gb}+ e^{v_n}}$ is bounded independently of $n$, it follows by the dominated convergence theorem that 
\begin{equation}\label{eq 4.8}
\int_{\BBR^2}W_\gb\frac{ e^{v_{\beta,A}}}{e^{\nu_1-\nu_2}\gl^{-\gb}+ e^{v_{\beta,A}}}
dx=\int_{\BBR^2}g_\gb dx.
\end{equation}
Combining this identity with the estimate
$$\left|g_\gb(x)-W_\gb(x)\frac{   e^{ v_{\beta,A}(x)}}{(e^{\nu_1-\nu_2}\lambda^{-\beta}+ e^{v_{\beta,A}(x)})^{1+a}}\right|
\leq c_{43}(1+|x|)^{-2aN-a\gb},
$$
and using Lemma \ref{lm 1.2}, we infer that $v_{\beta,A}$ is uniformly bounded in $\BBR^2$ and that there exists $C_{\gb, A}>A$ such that 
\begin{equation}\label{eq 4.9}
v_{\gb,A}(x)=C_{\gb, A}+O(|x|^{-\frac{a\gb+2aN-2}{a\gb+2aN-1}})\quad\text{as }\;|x|\to+\infty.
\end{equation}

In order to construct the sequence of solutions, we start with $A=A_0=1$, then take $A=A_1=\inf\{k\in\BBN:k>C_{\gb, A}\}$ and we iterate this process, defining by induction $A_{i+1}$ by
$A_{i+1}=\inf\{k\in\BBN:k>C_{\gb, A_i}\}$.

 \hfill$\Box$

\medskip

\noindent{\bf  Proof of Theorem \ref{teo 0} part $(ii)$ and Theorem \ref{teo 1}.}
 {\it Multiple solutions.} Let
 $$u_{\beta,i} =-\nu_1 +\nu_2 +\beta \ln \lambda +v_{\beta,i},$$
where $\{v_{\beta,i}\}_i$ is a sequence solutions of (\ref{eq 4.1}) which exist by Proposition \ref{pr 4.1}. Then $\{u_{\beta,i}\}_i$ is a sequence of   non-topological solutions in type II of (\ref{eq 1.1}) verifying   (\ref{1.3}) and with total magnetic flux   $2\pi[2(N-M)+\beta]$.
The proof is now complete. \hfill$\Box$

\subsection{Topological solution }

\noindent{\bf Proof of Theorem \ref{teo 2}.} {\it Multiple Topological solutions.} Let $u$ be a  topological  solution   of
problem (\ref{eq 1.1}). We can write it  as $ u =-\nu_1 +\nu_2 +v$ 
where $v $ is a bounded regular solution of
\begin{equation}\label{eq 5.1}
 \displaystyle  -\Delta v +V\frac{   e^{ v}}{(e^{\nu_1-\nu_2} + e^{v})^{1+a}} =g_0 \quad{\rm in}\;\; \mathbb{R}^2
\end{equation}
with
$$
 g_0 =f_1-f_2,
$$
and where the functions $f_1$ and $f_2$ have been defined in (\ref{V}). They are smooth, have compact support in $B_{r_0}(0)$ and the flux identity 
(\ref{Vf}) is satisfied.\smallskip

\noindent {\it Claim}: 
 Problem (\ref{eq 5.1}) possesses a sequence of bounded solutions  $\{v_{i}\}_i$ such that
  \begin{equation}\label{5.2}
  \int_{\R^2}V\frac{   e^{ v_{ i}}}{(e^{\nu_1-\nu_2} + e^{v_{i}})^{1+a}} dx =4\pi (N-M),
  \end{equation}
  and
\begin{equation}\label{e 5.3}v_{i}(x) = C_i+O(|x|^{-\frac{2aN-2}{2aN-1}})\quad {\rm as} \ \    |x|\to+\infty,
\end{equation}
 with $C_i\to+\infty$ as $ i\to+\infty.$\smallskip

\noindent  This can be proved as follows: given $A>0$,  let $w_A$ be the bounded solution  of
\begin{equation}\label{e 5.0}
-\Delta w+   e^{-A(1+a)} V    e^{ w}  =g_0\quad{\rm in}\ \,\R^2.
\end{equation}
We note that
$$w_A=w_0+A(1+a),$$
where $w_0$ is a bounded solution of (\ref{eq 4.2}) with $A=0$.
Note also that  if $A\geq A^*=a^{-1}\norm{w_0}_{L^\infty(\R^2)}$, then
$$w_A\geq A\quad{\rm in}\quad \R^2.$$
We set  $\mu_0=w_A$, and  define $\mu_n$ ($n\in\N$) to be the solution of
\begin{equation}\label{e 5.1}
-\Delta \mu_n+ V\frac{   e^{ \mu_n}}{(e^{\nu_1-\nu_2} + e^{\mu_{n-1}})^{1+a}} =g_0\quad{\rm in}\ \,\R^2.
\end{equation}
As in the proof of Proposition \ref{pr 4.1}, the mapping $n\mapsto \mu_n$ is increasing and $\mu_n$ is uniformly upper bounded. It converges to some function 
$v_A$ as $n\to+\infty$, and $v_A$ is a weak solution of (\ref{eq 5.1}). Since $V(x)\leq c_{22}|x|^{-2aN}$ when $|x|\to+\infty$, and $2aN>2$, there holds 
$$\mu_n(x)=C_{n,A}+O\left(|x|^{-\frac{2aN-2}{2aN-1}}\right)\quad\text{and }\;\left|\nabla \gm_n(x)\right|\leq c_{44}|x|^{-1-\frac{2aN-2}{2aN-1}}
\quad\text{as }\; |x|\to+\infty.
$$
Integrating (\ref{e 5.1}) on $B_R$ and letting $R\to\infty$ yields
$$\int_{\BBR^2}V\frac{   e^{ \mu_n}}{(e^{\nu_1-\nu_2} + e^{\mu_{n-1}})^{1+a}} dx=\int_{\BBR^2}(f_1-f_2)dx=4\gp(N-M). 
$$
Because $\frac{   e^{ \mu_n}}{(e^{\nu_1-\nu_2} + e^{\mu_{n-1}})^{1+a}}$ is uniformly bounded and $V\in L^1(\BBR^2)$ we obtain by the dominated convergence theorem
\begin{equation}\label{e 5.2}
\int_{\BBR^2}V\frac{   e^{ v_A}}{(e^{\nu_1-\nu_2} + e^{v_A})^{1+a}} dx=4\gp(N-M). 
\end{equation}
Therefore, 
\begin{equation}\label{e 5.3bis}
v_A(x)=C_A+O(|x|^{-\frac{2aN-2}{2aN-1}})\quad {\rm as} \ \    |x|\to+\infty,
\end{equation}
and the end of the proof is similar as the one of Proposition \ref{pr 4.1}.

\hfill$\Box$

\setcounter{equation}{0}
\section{  Nonexistence }

\begin{lemma}\label{lm 6.1}
Let  $aN< 1$. Then\smallskip

 \noindent $(i)$  Problem (\ref{eq 1.1}) has no  solution $u_\beta$ verifying
 \begin{equation}\label{6.3}
  u_\beta(x) -\beta\ln|x|=o(\ln|x|)\quad {\rm as} \ \    |x|\to+\infty
 \end{equation}
\qquad\ \  for $\beta^*<\beta\le 0$.\smallskip

 \noindent   $(ii)\ $   Problem (\ref{eq 1.1}) has no solution $u_\beta$  verifying (\ref{6.3})
 if $0\le \beta< \frac{2-aN}{a}$.\smallskip

  \noindent   $(iii)$    Problem (\ref{eq 1.1}) has no  topological  solution.
\end{lemma}
{\bf Proof.} We recall that a solution verifying (\ref{6.3}) with $\gb<0$ (resp. $\gb>0$) is called non-topological of type II (resp. type I). Given a function, we denote by $\overline w$  the circular average of $w$, i.e.
$$\overline w(r)=\frac1{2\pi r}\int_{\partial B_r(0)} w(\xi)d\theta(\xi)=\frac1{2\pi}\int_{0}^{2\gp} w(r,\gth)d\theta.$$
For $|x|\geq r_0$, there exists $c_{45}>0$ such that
$${\bf P}(x)\geq c_{45}|x|^{-aN},$$
and we set, for all $x\in\R^2$
\begin{equation}\label{6.1}
 h_u(x )={\bf P}(x)\frac{e^{u(x)}}{(1+e^{u(x)})^{1+a}}.
\end{equation}

\noindent  {\it Part $(i)$.}
If $u$ is a non-topological solution of Type I, it satisfies $u(x)\le c_{46}$  for $|x|\geq r_0$ and $c_{46}>0$. By Jensen inequality there exists positive constants $c_{47}$ and $c_{48}$ such that
\begin{equation}\label{6.4}
\overline h_u(r)\geq  \frac{ c_{47}A_0r^{-aN}}{(1+e^{c_{24}})^{1+a}}\overline{ e^{u(x)}} \geq c_{48} e^{\overline u(r)}\quad {\rm for}\ \ r> r_0,
\end{equation}
and from (\ref{6.3}),  there exist $\epsilon_0\in(0,\,1)$ and $c_{49}>0$ such that for $r> r_0$,
\begin{eqnarray}\label{6.2}
\overline h_u(r)  \geq      \frac{ c_{49} }{ r^{2-\epsilon_0} }.
\end{eqnarray}
Then (\ref{eq 1.1}) implies that
$$
  (r\overline u_r)_r  \geq   \frac{ c_{49} }{ r^{1-\epsilon_0} } \quad {\rm for}\ \, r\geq r_0,
$$
thus, integrating the above inequalities, we obtain 
\begin{eqnarray*}
 r\overline u_r(r)- r_0\overline u_r(r_0)  \geq    c_{50} (r^{\epsilon_1}-r_0^{\epsilon_1}),
\end{eqnarray*}
where $\overline u_r=\frac{d\overline u}{dr}$ and $c_{49},\,c_{50}>0$. Hence  there holds
$$
\overline u(r) \geq \overline u(r_0)+(r_0\overline u_r(r_0)-c_{50}r_0^{\epsilon_1}) \ln r+\frac{c_{50}}{\epsilon_1}r^{\epsilon_1} \quad{\rm for} \  r>r_0.
$$
As a consequence,
\begin{equation}
\overline u(r)  \to+\infty \quad{\rm as}\quad r\to +\infty,
\end{equation}
which contradicts the fact that  $u$ is bounded from above.\smallskip

\noindent  {\it Part $(ii)$.}
If $u$ is a non-topological solution of Type II and $0\le \beta< \frac{2-aN}{a}$,   then
\begin{eqnarray*}
 h_u(x ) &\geq & {\bf P}(x)\frac{1}{(1+e^u)^{a}}
   \geq \frac{ c_{51}  }{ |x|^{2-\epsilon_1} }
\end{eqnarray*}
for some $\epsilon_1>0$ and $c_{51}>0$.
Then (\ref{eq 1.1}) implies that
$$
  (r\overline w_r)_r  \geq   \frac{ c_{52} }{ r^{1-\epsilon_1} } \quad {\rm for}\ \ r\geq r_0.
$$
Hence  there holds
$$
\overline w(r) \geq \overline w(r_0)+(r_0\overline w_r(r_0)-c_{52}r_0^{\epsilon_1}) \ln r+\frac{c_{52}}{\epsilon_1}r^{\epsilon_1} \quad{\rm for} \ \, r>r_0,
$$
which contradicts (\ref{6.3}). \smallskip

\noindent  {\it Part $(iii)$.} The proof is the same as above.  \hfill$\Box$\medskip

\noindent{\bf Proof of Theorem \ref{teo 3}.} If $aN<1$, Lemma \ref{lm 6.1} implies that
then   problem (\ref{eq 1.1}) has no solution $u_\beta$ for $\beta^*<\beta< \frac{2-aN}{a}$ verifying  $u_\beta(x)=\beta\ln|x|+O(1)$.   \smallskip

  Next we assume that   $aN= 1$, and  $u$ is a  topological solution  (\ref{eq 1.1}).
Hence $u$ is bounded at infinity and
\begin{eqnarray*}
 h_u(x ) &\geq & {\bf P}(x)\frac{e^u}{(1+e^u)^{a}}
   \geq \frac{ c_{53}  }{ |x|^{-2 }}.
\end{eqnarray*}
Then (\ref{eq 1.1}) implies that
$$
  (r\overline u_r)_r  \geq   \frac{ c_{54} }{ r  } \quad {\rm for}\ \, r\geq r_0.
$$
By integrating this inequality we encounter a contradiction with the fact that $u$ is bounded at infinity.
\hfill$\Box$

\bigskip

 \noindent{\small {\bf Acknowledgements:}  The authors are grateful to the referees for careful checking of the manuscript and relevant observations concerning some 
theoritical aspects of the physical theory. \\
H. Chen  is supported by NSFC (No:11726614, 11661045).

\end{document}